\documentclass[12pt,a4paper,reqno]{amsart}
\usepackage{amsmath}
\usepackage{amsfonts}
\usepackage{amssymb}
\numberwithin{equation}{section}

\theoremstyle{plain}

\newtheorem{theorem}[subsection]{Theorem}
\newtheorem{proposition}[subsection]{Proposition}
\newtheorem{lemma}[subsection]{Lemma}
\newtheorem{corollary}[subsection]{Corollary}
\newtheorem{conjecture}[subsection]{Conjecture}

\theoremstyle{definition}

\newtheorem{definition}[subsection]{Definition}

\newtheorem{question}[subsection]{Question}

\renewcommand{\leq}{\leqslant}
\renewcommand{\geq}{\geqslant}

\newsavebox{\proofbox}
\savebox{\proofbox}{\begin{picture}(7,7)%
  \put(0,0){\framebox(7,7){}}\end{picture}}
\def\boxeq{\tag*{\usebox{\proofbox}}}

\newcommand{\md}[1]{\ensuremath{(\mbox{mod}\, #1)}}

\def\proof{\noindent\textit{Proof. }}
\def\endproof{\hfill{\usebox{\proofbox}}}

\def\Agh{A_H^{+g}}
\def\ind{\mbox{ind}}

\def\Agh{A_H^{+g}}

\def\SBA{\widetilde{B}_{\Gamma,\delta}}

\def\sba{\beta_{\Gamma,\delta}}

\def\psa{\psi_{\Gamma,\delta}}

\def\psb{\psi_{\Gamma',\delta'}}

\begin{document}

\title[Szemer\'edi-type regularity lemma]{A Szemer\'edi-type regularity lemma in abelian groups, with applications
}

\author{Ben Green}
\address{
}
\email{bjg23@hermes.cam.ac.uk}

\thanks{While this work was carried out, the author was supported by a fellowship of Trinity College, Cambridge, England and a PIMS postdoctoral fellowship at the University of British Columbia, Vancouver, Canada}

\begin{abstract}
Szemer\'edi's regularity lemma is an important tool in graph theory which has applications throughout combinatorics. \vspace{8pt}

\noindent In this paper we prove an analogue of Szemer\'edi's regularity lemma in the context of abelian groups and use it to derive some results in additive number theory.\vspace{8pt}

\noindent One is a structure theorm for sets which are almost sum-free. If $A \subseteq \{1,\dots,N\}$ has $\delta N^2$ triples $(a_1,a_2,a_3)$ for which $a_1 + a_2 = a_3$ then $A = B \cup C$, where $B$ is sum-free and $|C| = \delta' N$, and $\delta' \rightarrow 0$ as $\delta \rightarrow 0$\vspace{8pt}

\noindent Another answers a question of Bergelson, Host and Kra. If $\alpha, \epsilon > 0$, if $N > N_0(\alpha,\epsilon)$ and if $A \subseteq \{1,\dots,N\}$ has size $\alpha N$, then there is some $d \neq 0$ such that $A$ contains at least $(\alpha^3 - \epsilon)N$ three-term arithmetic progressions with common difference $d$.

\end{abstract}

\maketitle

\section{Introduction}\label{sec1} \noindent Szemer\'edi's regularity lemma \cite{Szem} is an important result in graph theory with numerous applications in combinatorics and number theory. It has been described as a structure theorem for an arbitrary graph. We give a very brief introduction to the regularity lemma which is designed to motivate the results of this paper. The reader may find a much more extensive survey in the excellent article of Koml\'os and Simonovits \cite{KomSim}, and a nicely-explained proof of the lemma in \cite{Bol}.\\[11pt]
Let $\Gamma = (V,E)$ be a graph and let $A,B$ be disjoint subsets of $V$. Define the density $d(A,B)$ to be the proportion of elements $(x,y) \in A \times B$ such that $xy \in E(\Gamma)$. If $\epsilon > 0$, we say that a pair $(A,B)$ is $\epsilon$-uniform if 
\[ |d(A',B') - d(A,B)| \; \leq \; \epsilon \] whenever $A' \subseteq A$ and $B' \subseteq B$ satisfy $|A'| \geq \epsilon |A|$ and $|B'| \geq \epsilon |B|$.
\begin{proposition}[Szemer\'edi's regularity lemma]\label{prop1.1} Let $\epsilon > 0$. There exists $M = M(\epsilon)$ such that the vertex set $V(\Gamma)$ of any graph $\Gamma$ can be partitioned into $\lceil 1/\epsilon \rceil \leq m \leq M$ sets $V_1,\dots,V_m$ with sizes differing by at most 1, such that at least $(1 - \epsilon)m^2$ of the pairs $(V_i,V_j)$ are $\epsilon$-uniform.
\end{proposition} 
\noindent Hereafter we will refer to Szemer\'edi's regularity lemma as SzRL.
One reason that SzRL has been described as a ``structure theorem for all graphs'' is the fact that it is possible to say much more about the bipartite graph induced by a regular pair $(V_i,V_j)$ than it is about an arbitrary graph. As an example of this phenomenon we cite the following result.
\begin{proposition}[Counting lemma] \label{prop1.2}Suppose that $U,V$ and $W$ are disjoint sets of $s$ vertices in some graph $\Gamma$. Write $d(U,V) = \gamma$, $d(V,W) = \alpha$ and $d(U,W) = \beta$, and suppose that the pairs $(U,V)$ and $(V,W)$ are $\epsilon$-regular. Write $T(U,V,W)$ for the number of triangles in $U \times V \times W$ \emph{(}that is, triples $(u,v,w)$ such that $uv,vw,wu \in E(\Gamma)$\emph{)}. Then 
\[ \left| T(U,V,W) - \alpha\beta\gamma s^3 \right| \; \leq \; 4\epsilon s^3.\]
\end{proposition}
\noindent Observe that only \textit{two} of the pairs $(U,V)$, $(V,W)$, $(U,W)$ are required to be regular. Proposition \ref{prop1.2} may be combined with SzRL to prove the following.
\begin{proposition} \label{prop1.3}Let $\Gamma$ be a graph on $n$ vertices, and suppose that $\Gamma$ contains $o(n^3)$ triangles. Then we may remove $o(n^2)$ edges from $\Gamma$ so as to leave a graph which is triangle-free.
\end{proposition}
\noindent The notation here is convenient but offers scope for confusion. What we mean is that there is a function $\delta' = \delta'(\delta)$ such that $\delta' \rightarrow 0$ as $\delta \rightarrow 0$, and which has the following property. If $\Gamma$ contains at most $\delta n^3$ triangles then we may remove $\delta' n^2$ edges from $\Gamma$ so as to leave a graph which is triangle-free. \\[11pt]
We have not attributed Proposition \ref{prop1.3}, as it is not clear to us where it was first stated. A slightly weaker result was obtained by Ruzsa and Szemer\'edi in 1976 \cite{RuzSzem}. At that time SzRL had only been formulated for bipartite graphs. The modern \cite{Szem} formulation, together with the ideas of \cite{RuzSzem}, would certainly imply Proposition \ref{prop1.3}. The result is also well-known in the literature concerning ``property testing'': see, for example, \cite{Alo}.\\[11pt]
Proposition \ref{prop1.3} is surprising and interesting in its own right. It also has important applications, not the least of which is a simple proof of Roth's theorem that $r_3(n)$, the size of the largest subset of $\{1,\dots,n\}$ containing no 3-term arithmetic progression, satisfies $r_3(n) = o(n)$.\\[11pt]
One of the results of this paper is an analogous theorem for abelian groups. Let $G$ be an abelian group with cardinality $N$, and let $A \subseteq G$. A triple $(x,y,z) \in A^3$ is a  \textit{triangle} if $x + y + z = 0$. 
\begin{theorem}\label{thm1.4}
Suppose that $A \subseteq G$ is a set with $o(N^2)$ triangles. Then we may remove $o(N)$ elements from $A$ to leave a set which is triangle-free.\end{theorem}
\noindent In fact, we will deduce this result from the following more general theorem.
\begin{theorem}\label{thm1.5}
Let $k \geq 3$ be a fixed integer, and suppose that $A_1,\dots,A_k$ are subsets of $G$ such that there are $o(N^{k-1})$ solutions to the equation $a_1 + \dots + a_k = 0$ with $a_i \in A_i$ for all $i$. Then we may remove $o(N)$ elements from each $A_i$ so as to leave sets $A'_i$, such that there are no solutions to $a'_1 + \dots + a'_k = 0$ with $a'_i \in A'_i$ for all $i$.
\end{theorem}
\noindent A simple corollary is the structure theorem for sets of integers which are almost sum-free, as featured in the abstract of the paper.
\begin{corollary}\label{cor1.6}
Suppose that $A \subseteq [N]$ is a set containing $o(N^2)$ triples with $x + y = z$. Then $A = B \cup C$ where $B$ is sum-free and $|C| = o(N)$.
\end{corollary}
\noindent The proof of Theorems \ref{thm1.4} and \ref{thm1.5} are in many ways analagous to the proof of Proposition \ref{prop1.3}. In particular we must prove a regularity lemma in the context of abelian groups. Although this regularity lemma (Theorem \ref{thm5.3}) is probably the most interesting result in the paper, it takes some time to set up the notation necessary to state it and so we do not do so here. We will, however, give a sketch of how Proposition \ref{prop1.3} follows from SzRL and the counting lemma (Proposition \ref{prop1.2}). The deduction of Theorem \ref{thm1.4} from Theorem \ref{thm1.5} and an appropriate analogue of the counting lemma is in many ways quite similar. We start with a definition.
\begin{definition}[Reduced graph]\label{def1.7} Let $\Gamma = (V,E)$ be a graph and let $\epsilon > 0$. Take a partition $V_1 \cup \dots \cup V_m$ satisfying the conclusions of SzRL. Define a new graph $\Gamma'$ as follows. For each pair $(i,j)$, consider whether one of the following three conditions is satisfied:
\begin{enumerate}
\item $i = j$;
\item $(V_i,V_j)$ is not $\epsilon$-regular;
\item $d(V_i,V_j) < 2\epsilon^{1/3}$.
\end{enumerate}
If so, delete all edges from $V_i$ to $V_j$. Let $\Gamma'$ be the graph that remains; we refer to $\Gamma'$ as an $\epsilon$-reduced subgraph of $\Gamma$. Often, when the parameters are clear from the context, we will refer to $\Gamma'$ simply as the \emph{reduced graph}.
\end{definition}
\noindent It is not hard to see that if $\Gamma$ has $N$ vertices then
\begin{equation}\label{eq1} |E(\Gamma')| \; \geq \; |E(\Gamma)| - 10\epsilon^{1/3} N^2.\end{equation}
\textit{Sketch proof of Proposition \ref{prop1.3}.} Define, for each $\delta \in (0,1]$, a value $\epsilon = \epsilon(\delta)$ for which $4\epsilon(M(\epsilon))^{-3} > \delta$, but so that $\epsilon(\delta) \rightarrow 0$ as $\delta \rightarrow 0$. Suppose that $\Gamma$ is a graph with $\delta n^3$ triangles. Consider $\Gamma'$, an $\epsilon$-reduced subgraph of $\Gamma$ relative to some underlying partition $V_1 \cup \dots \cup V_m$ coming from SzRL. We know from \eqref{eq1} that $\Gamma'$ is obtained from $\Gamma$ by the deletion of relatively few edges. We claim that $\Gamma'$ is triangle-free. If this is not the case then it contains a triangle $v_iv_jv_k$ with $v_i \in V_i$, $v_j \in V_j$, $v_k \in V_k$. For simplicity assume that $|V_i| = |V_j| = |V_k|$. Each of these sets has size at least $n/2M(\epsilon)$. Now by the construction of $\Gamma'$ we see that $(V_i,V_j)$ is $\epsilon$-regular and $d(V_i,V_j) \geq 2\epsilon^{1/3}$, and similarly for $(V_j,V_k)$ and $(V_k,V_i)$. Thus, by the counting lemma, we see that 
\[ T(V_i,V_j,V_k) \; \geq \; 4\epsilon(M(\epsilon))^{-3}n^3 \; > \; \delta n^3,\] contrary to assumption.\endproof\\[11pt]
Our proof of Theorem \ref{thm1.4} will follow a broadly similar scheme. There will be a regularity lemma, a counting lemma and, given a set $A \subseteq G$, a definition of an $\epsilon$-reduced sub\textit{set} $A'$ of $A$. \\[11pt]
A somewhat different application of our regularity lemma is to a question raised by Bergelson, Host and Kra in \cite{BHK}. The question was this:
\begin{question}\label{BHK-question} Suppose that $\alpha,\epsilon > 0$. Is it true that there is $N_0(\alpha,\epsilon)$ such that if $N > N_0(\alpha,\epsilon)$, and if $A \subseteq \{1,\dots,N\}$ has size $\alpha N$, then there is some $d \neq 0$ such that $A$ has at least $(\alpha^3 - \epsilon)N$ three-term arithmetic progressions with common difference $d$?
\end{question}
\noindent In fact the question was also asked for 4-term progressions, and it was shown that no such result holds for progressions of length 5 and higher. We do not answer the question about 4-term progressions here (though see the remarks in \S \ref{sec8} IV). \vspace{11pt}

\noindent Question \ref{BHK-question} is addressed in \S \ref{sec7.5}. There we begin with the following result, whose proof is a relatively clean application of our regularity lemma.

\begin{theorem}\label{thm1.8} 
Suppose that $\alpha,\epsilon > 0$. Then there is $N_0(\alpha,\epsilon)$ such that if $G$ is an abelian group of size $N > N_0(\alpha,\epsilon)$ with $N$ odd, and if $A \subseteq G$ has size $\alpha N$, then there is some $d \neq 0$ such that $A$ has at least $(\alpha^3 - \epsilon)N$ three-term arithmetic progressions with common difference $d$.
\end{theorem}
\noindent This does not seem, as it stands, to give an affirmative answer to Question \ref{BHK-question}. By modifying the argument in some small but slightly technical ways, we can answer that question.
\begin{theorem}\label{thm1.9} The answer to Question \ref{BHK-question} is yes.\end{theorem}
\noindent Let us conclude this introduction with a word or two on notation. Let $G$ be a finite abelian group and let $G^{\ast}$ be the dual of $G$, thought of as the group of characters $\gamma : G \rightarrow \mathbb{C}$. If $f : G \rightarrow \mathbb{R}$ is a function and $\gamma \in G^{\ast}$ a character, define the Fourier transform $\widehat{f}(\gamma) = \sum_x f(x)\gamma(x)$. Sometimes, when taking the Fourier transform of a reasonably complicated expression, we will use the alternative notation $(\mbox{expression})^{\wedge}(\gamma)$. If $f,g : G \rightarrow \mathbb{R}$ define the convolution $f \ast g(x) = \sum_y f(y)g(x - y)$. A number of simple instances of Young's inequality, such as the bounds $\Vert f \ast g \Vert_1 \leq \Vert f \Vert_1 \Vert g \Vert_1$ and $\Vert f \ast g \Vert_{\infty} \leq \Vert f \Vert_{\infty} \Vert g \Vert_1$, will be used without comment.

\section{A study of the group $(\mathbb{Z}/2\mathbb{Z})^n$}\label{sec2} 
\noindent It is rather hard to describe the regularity lemma for a general group $G$ (that will be the objective of sections \ref{sec3}, \ref{sec4} and \ref{sec5}). The group $\mathbb{Z}/2N\mathbb{Z}$, which is of interest as regards, say, Corollary \ref{cor1.6}, has all the difficulties of the general case. It turns out, however, that everything works out very cleanly in the particular case $G = (\mathbb{Z}/2\mathbb{Z})^n$. This is one more instance (cf. \cite{Gre1,Mes,Ruz2}) in which the consideration of vector spaces over finite fields facilitates thinking about questions concerning the integers. For a survey of this phenomenon, the reader may care to consult the survey article \cite{green-survey}.\\[11pt]
This section is devoted to this special case, and is independent of the rest of the paper. What we describe is possibly the very simplest situation in which regularity and associated ideas such as the counting lemma can be studied.\\[11pt]
For the remainder of \S \ref{sec2} set $G = (\mathbb{Z}/2\mathbb{Z})^n$ and write $N = |G| = 2^n$. 
Let $H \leq G$ be a subgroup. For any $g \in G$ we may define a set $\Agh \subseteq H$ by setting
\[ \Agh(x) \; = \; A(x + g)\] for $x \in H$. These sets represent intersections of $A$ with cosets of $H$. We will be interested in the Fourier coefficients of $\Agh$, defined for $\eta \in H^{\ast}$ by
\[ \widehat{\Agh}(\eta) \; = \; \sum_{x \in H} \Agh(x) (-1)^{\langle x,\eta\rangle}.\] If 
\[ \sup_{\eta \neq 0} |\widehat{\Agh}(\eta)| \; \leq \; \epsilon |H|\] then we say that $g$ is an $\epsilon$-regular value with respect to $A$ (and the subgroup $H$). \\[11pt]
The concept of regularity in this sense has been well-studied as a notion of pseudorandomness for subsets of abelian groups. See \cite{ChuGra,Gow2} for more details, other applications and equivalent formulations.\\[11pt]
If the number of $g \in G$ which fail to be $\epsilon$-regular is no more than $\epsilon N$ then we say that the subgroup $H$ is $\epsilon$-regular for $A$. In the statement of the following result, $W(t)$ is defined to be a tower of twos of height $\lceil t\rceil$.
\begin{theorem}[Regularity lemma in $(\mathbb{Z}/2\mathbb{Z})^n$]\label{thm2.1}
Let $\epsilon \in (0,\frac{1}{2})$ and let $A \subseteq G$. Then there is a subgroup $H \leq G$ of index at most $W(\epsilon^{-3})$ which is $\epsilon$-regular for $A$.
\end{theorem}
\proof  If $H \leq G$, and if the sets $\Agh$ are as above, define a quantity $\ind(A;H)$ by
\[ \ind(A;H) \; = \; \frac{1}{N}\sum_{g \in G} \left(\frac{|\Agh|}{|H|}\right)^2.\]Observe that $0 \leq \ind(A;H) \leq 1$ for any $A$ and $H$. 
We will define a sequence of subgroups
\[ G \; = \; H_0 \; \geq \; H_1 \; \geq \; H_2 \; \geq \; H_k\] of (very rapidly) increasing index. If $H_i$ is not $\epsilon$-regular for $A$ then we will use any abnormally large Fourier coefficients of the sets $A_{H_i}^{+g}$ to construct an $H_{i+1}$ for which $\ind(A,H_{i+1})$ is substantially larger than $\ind(A,H_i)$. This process must terminate after a finite number of steps (depending on $\epsilon$).\\[11pt]
The heart of the proof is the process of passing from $H_i$ to $H_{i+1}$. This is covered in detail in the following lemma.
\begin{lemma}\label{lem2.2} Let $\epsilon \in (0,\frac{1}{2})$ and
suppose that $H \leq G$ is a subgroup which is not $\epsilon$-regular for $A$. Then there is a subgroup $H' \leq H$ such that $|G/H'| \leq 2^{|G/H|}$ and $\mbox{\emph{ind}}(A,H') \geq \mbox{\emph{ind}}(A,H) + \epsilon^3$.
\end{lemma}
\proof There are at least $\epsilon N$ values of $g \in G$ for which $\sup_{\eta \neq 0} |\widehat{\Agh}(\eta)| \geq \epsilon |H|$. Now if $g_1$ and $g_2$ lie in the same coset of $H$ then the sets $A_H^{+g_1}$ and $A_H^{+g_2}$ are just translates of one another. In particular it is easy to see that 
\[ \widehat{A_H^{+g_1}}(\eta) \; = \; (-1)^{\langle g_2 - g_1 ,\eta\rangle} \widehat{A_H^{+g_2}}(\eta),\] and so $A_H^{+g_1}$ and $A_H^{+g_2}$ have large Fourier coefficients at the same points. This means that there are cosets $H + g_i$, $i = 1,\dots,K$, $\frac{1}{2}|G/H| \geq K \geq \epsilon|G/H|$ and points $\eta_1,\dots,\eta_K \in H^{\ast}$ such that $|\widehat{\Agh}(\eta_i)| \geq \epsilon |H|$ for all $g \in H + g_i$.\\[11pt]
Let $H' \leq H$ be the annihilator of the $\eta_i$, that is to say the set of  $x$ such that $\langle x,\eta_i\rangle = 0$ for all $i = 1,\dots,K$. The bound on $|G/H'|$ claimed in the lemma is immediate, and we must check that $\ind(A,H') \geq \ind(A,H) + \epsilon^3$ as stated. As a first observation, note that $\widehat{H'}(0)$ and all of the $\widehat{H'}(\eta_i)$ are equal to $|H'|$.\\[11pt]
Now one has
\begin{eqnarray*}
N|H'|^2|H| \cdot \ind(A,H') \; = \; \sum_g |A^{+g}_{H'}|^2  & = &  \sum_{g \in G} \sum_{h \in H} |A^{+(g+h)}_{H'}|^2 \\ \nonumber & = &  \sum_{g \in G} \sum_{h \in H}\left|\sum_{x \in G} A(x - g - h)H'(x)\right|^2 \\ 
& = & \sum_{g \in G} \sum_{h \in H} (\Agh \ast H')(h)^2.
\end{eqnarray*}
The sum over $h \in H$ may be written in terms of the Fourier transform on the subgroup $H$. One has
\[
 N|H|^2|H'|^2 \ind(A,H') \; = \; \sum_{g \in G}\sum_{\eta \in H^{\ast}}|\widehat{\Agh}(\eta)|^2|\widehat{H'}(\eta)|^2.\]
The term with $\eta = 0$ is easily seen to be $N|H|^2|H'|^2 \ind(A,H)$. To bound the sum over $\eta \neq 0$ from below, write
\begin{eqnarray*}
\sum_{g \in G}\sum_{\eta \neq 0}|\widehat{\Agh}(\eta)|^2|\widehat{H'}(\eta)|^2 & \geq & \nonumber \sum_{i =1}^K \sum_{g \in H+ g_i} |\widehat{\Agh}(\eta_i)|^2|\widehat{H'}(\eta_i)|^2 \\ \nonumber & \geq & \epsilon^2 K |H|^3|H'|^2 \\ & \geq & \epsilon^3 N|H|^2|H'|^2.\end{eqnarray*}
These observations confirm that $\ind(A,H') \geq \ind(A,H) + \epsilon^3$, which is what we set out to prove.\endproof\\[11pt]
To deduce Theorem \ref{thm2.1}, simply carry out the program outlined just before the statement of Lemma \ref{lem2.2}. Set $H_0 = H$, and define subgroups $H_i$ inductively. If $H_i$ is not $\epsilon$-regular for $A$, apply Lemma \ref{lem2.2} with $H = H_i$, and set $H_{i+1} = H'$. It is clear that $\ind(A,H_i) \geq i\epsilon^3$, and so this process can take place no more than $\lfloor \epsilon^{-3} \rfloor$ times. In that time, the index $|G/H_i|$ has not become more than $W(\epsilon^{-3})$.\\[11pt]
The reader who is familiar with the proof of SzRL will notice some strong similarities between that argument and the proof of Theorem \ref{thm2.1}.\\[11pt]
The next result is a counting lemma directly analagous to Proposition \ref{prop1.2}.
\begin{proposition}[Counting lemma in $(\mathbb{Z}/2\mathbb{Z})^n$.] \label{prop2.3}
Suppose that $H$ is a subgroup of $G$, and let $g_1,g_2,g_3 \in G$. Suppose that $|A^{+g_i}_H| = \alpha_i |H|$, and suppose that $A^{+g_1}_H$ is $\epsilon$-regular. Then $T(g_1,g_2,g_3)$, the number of triples $(x_1,x_2,x_3)$ such that $x_i \in A^{+g_i}_H$ and $x_1 + x_2 + x_3 = 0$, satisfies 
\[ \left|T(g_1,g_2,g_3) - \alpha_1\alpha_2\alpha_3|H|^2\right| \; \leq \; \epsilon|H|^2.\]
\end{proposition}
\proof One has, using orthogonality relations for characters,
\begin{eqnarray*} T(g_1,g_2,g_3) & = & \sum_{\substack{ x_1,x_2,x_3 \in H \\ x_1 + x_2 + x_3 = 0}} A^{+g_1}_H(x_1)A^{+g_2}_H(x_2)A^{+g_3}_H(x_3).\\ & = & |H|^{-1}\sum_{\gamma} \widehat{A^{+g_1}_H}(\gamma)\widehat{A^{+g_2}_H}(\gamma)\widehat{A^{+g_3}_H}(\gamma).\end{eqnarray*}
The term with $\gamma = 0$ is precisely $\alpha_1\alpha_2\alpha_3|H|^2$. One can estimate the remainder with a simple $\ell^{2}$-$\ell^{\infty}$ inequality. Indeed
\begin{eqnarray*}
\left|\sum_{\gamma \neq 0} \widehat{A^{+g_1}_H}(\gamma)\widehat{A^{+g_2}_H}(\gamma)\widehat{A^{+g_3}_H}(\gamma)\right| & \leq & \sup_{\gamma \neq 0} \left|\widehat{A^{+g_1}_H}(\gamma)\right|\left(\sum_{\gamma}  \left|\widehat{A^{+g_2}_H}(\gamma)\right|^2\right)^{1/2}\left(\sum_{\gamma}  \left|\widehat{A^{+g_3}_H}(\gamma)\right|^2\right)^{1/2}.\\ & \leq & \epsilon |H|^2,\end{eqnarray*} the latter step following from the $\epsilon$-regularity of $A^{+g_1}_H$ and Parseval's identity. The result follows.\endproof\\[11pt]
Now we define a notion of reduced set, corresponding to the definition of reduced graph (Definition \ref{def1.7}).
\begin{definition}[Reduced set]\label{def2.4} Suppose that $A \subseteq G$, and let $H$ be $\epsilon$-regular for $A$. Define a new set $A'$ as follows. For each $g \in A$, consider whether either of the following two conditions is satisfied:\begin{enumerate}
\item $\Agh$ is not $\epsilon$-regular;
\item $|\Agh| \leq (2\epsilon)^{1/3}|H|$.
\end{enumerate}
If so, delete all of $\Agh$. Let $A'$ be the set that remains.
\end{definition}
\noindent Observe that (ii) depends only on the coset of $H$ that $g$ lies in. Write $X$ for the set of all $g$ satisfying (ii). Let $S \subseteq G/H$ be the set of cosets met by $X$, and for each $s \in S$ select some $g_s \in S \cap X$. Since $X = S \cap A$ we have
\[ |X| \; \leq \; \sum_{s \in S}|A^{+g_s}_H| \; \leq \; (2\epsilon)^{1/3}|S||H| \; \leq \; (2\epsilon)^{1/3}|G/H||H| \; = \; (2\epsilon)^{1/3}N.\]
The number of $g$ satisfying (i) is at most $\epsilon N$, and so
\begin{equation}\label{eq6} |A'| \; \geq \; |A| - 3\epsilon^{1/3}N.\end{equation}
\noindent We are now in a position to prove Theorem \ref{thm1.4} for $G$. Recall that a \textit{triangle} in a set $A$ is a triple $(x,y,z) \in A^3$ with $x + y + z = 0$.
\begin{theorem} Suppose that $A \subseteq (\mathbb{Z}/2\mathbb{Z})^n$ contains $o(N^2)$ triangles. Then we may remove $o(N)$ elements from $A$ to leave a set which is triangle-free.
\end{theorem}
\proof Define, for each $\delta \in (0,1]$, a value $\epsilon = \epsilon(\delta)$ for which $\epsilon (\Delta(\epsilon))^{-1/2} > \delta$, but so that $\epsilon(\delta) \rightarrow 0$ as $\delta \rightarrow 0$. Suppose that $A \subseteq G$ is a set with $\delta n^2$ triangles. Consider $A'$, an $\epsilon$-reduced subset of $A$ relative to some subgroup $H \leq G$ of index at most $\Delta(\epsilon)$ and which is $\epsilon$-regular for $A$. The existence of such an $H$ is the content of Theorem \ref{thm2.1}. We know from \eqref{eq6} that $A'$ is obtained from $A$ by the deletion of at most $3\epsilon^{1/3}N$ elements. We claim that $A'$ is triangle-free. Suppose that it contains three elements $y_1, y_2,y_3$ with $y_1 + y_2 + y_3 = 0$. But every triple $(x_1,x_2,x_3) \in H^3$ with $x_i \in A^{+y_i}_H$ and $x_1 + x_2 + x_3 = 0$ gives rise to a triangle $(x_1 + y_1,x_2 + y_2 ,x_3 + y_3)$ in $A$. By Proposition \ref{prop2.3} and the fact that the $x_i$ satisfy neither condition (i) nor (ii) in Definition \ref{def2.4}, the number of such triples is at least $\epsilon |H|^2$, which is more than $\delta N^2$.
This is contrary to assumption.\endproof\\[11pt]
In \S \ref{sec9} we will give an example in the spirit of Gowers \cite{Gow1} which shows that the huge bound which occurs in Theorem \ref{thm2.1} is to some extent necessary. For now, however, we press on with the main aim of the paper, which is the generalisation of the above to an arbitrary finite abelian group.

\section{The treatment of arbitrary abelian groups - introduction}\label{sec3} \noindent There is one rather obvious obstacle to generalising the results of \S \ref{sec2} to a general abelian group $G$: the lack, in general, of a plentiful supply of subspaces. In place of them, we will use \textit{Bohr neighbourhoods}.\\[11pt]
For the purposes of this paper, we will define the argument $\arg z$ of a complex number $z$ to lie in the interval $(-\pi,\pi]$. Let $\Gamma = \{\gamma_1,\dots,\gamma_d\}$ be a set of characters on $G$. Define the Bohr neighbourhood $B_{\Gamma,\delta}$ by
\[ B_{\Gamma,\delta} \; = \; \left\{x \; : \; |\arg \gamma_j(x)| \leq 2\pi\delta \;\; \mbox{for all $j = 1,\dots,d$}\right\}.\]
It is convenient to write $\Vert x \Vert_{\Gamma} = \sup_{j} |\frac{1}{2\pi}\arg \gamma_j(x)|$, so that $B_{\Gamma,\delta}$ is simply the set $\{x: \Vert x \Vert_{\Gamma} \leq \delta\}$.
When $G = (\mathbb{Z}/2\mathbb{Z})^n$, Bohr neighbourhoods are just subgroups. For other groups this is not the case, and in general there are some fairly substantial differences between their behaviour and that of a true subgroup. In $\mathbb{Z}/N\mathbb{Z}$ the Bohr neighbourhood $B = B_{\Gamma,\delta}$ tends to resemble a $d$-dimensional convex body, so that typically $|B + B|$ will be of cardinality closer to $2^d|B|$ than to $|B|$. Thus $B$ fails to behave like a group to a substantial extent, which (it turns out) means that there is no sensible way to do harmonic analysis on $B$, even approximately, in the manner we described in \S \ref{sec2}.\vspace{11pt}

\noindent The method we use to get around this stems from a beautiful observation of Bourgain \cite{Bou}. If $B' = B_{\Gamma,\delta'}$, where $\delta' \ll \delta$, then $B'$ tends to resemble a scaled-down version of $B$. It might then be expected (perhaps by thinking geometrically, imagining $B$ and $B'$ to be convex bodies) that $|B + B'| \approx |B|$. Roughly speaking \textit{pairs} of Bohr neighbourhoods, one much smaller than the other, are an appropriate substitute for subspaces. We will see in later sections that it is possible to do a sort of approximate harmonic analysis on such pairs of Bohr sets. 

\section{Properties of smoothed Bohr neighbourhoods}\label{sec4} \noindent In this section we define what may be called smoothed Bohr neighbourhoods and establish the basic properties of these functions that we will need. These are needed because there was one respect in which the discussion of \S \ref{sec3} was too simplistic. It turns out that $B_{\Gamma,(1 - \kappa)\delta}$ and $B_{\Gamma,\delta}$ can be quite different, even for very small $\kappa$. For a simple example, take $G = (\mathbb{Z}/5\mathbb{Z})^n$ and $\delta = 2\pi/5$. This kind of behaviour means that Bohr neighbourhoods do not always behave in a similar manner to convex bodies. \\[11pt]
Bourgain circumvented this obstacle by showing that for a fixed $\Gamma$, most values of $\delta$ are such that $B_{\Gamma,\delta}$ behaves in what he calls a \textit{regular} fashion. This makes the details of the argument even more difficult. In an exposition of Bourgain's work, Tao \cite{Tao} effects a significant simplification by putting this averaging over $\delta$ into the definition, getting a kind of smoothed Bohr neighbourhood. We give a different construction which is nonetheless inspired by this idea of Tao.\\[11pt]
This is a technical section of the paper the reader will lose little by simply looking at the definition of the functions $\psa$ (Definition \ref{def4.5}) and very briefly checking out the statements of their properties as laid down in Lemma \ref{lem4.5}. 
We begin with some simple properties of (unsmoothed) Bohr neighbourhoods.
\begin{lemma}\label{lem4.1} Let $G$ be an abelian group of size $N$, let $\Gamma = \{\gamma_1,\dots,\gamma_d\}$ be a set of $d$ characters on $G$ and let $\delta > 0$. Then\\
\emph{(i)} $|B_{\Gamma,\delta}| \geq \delta^d N$.\\
\emph{(ii)} $|B_{\Gamma,2 \delta}| \leq 5^d|B_{\Gamma,\delta}|$.
\end{lemma}
\proof If $t \in \mathbb{R}/\mathbb{Z}$ write $|t|$ for that representative of $t\pmod{1}$ which lies in the interval $(-1/2,1/2]$. For any $\eta$, write $S_{\eta}$ for the set of all for the set of all $y = (y_1,\dots,y_d) \in \mathbb{R}^d/\mathbb{Z}^d$ for which $|y_j| \leq \eta/2$ for all $j = 1,\dots,d$. Now if $x \in G$ write \[ v(x) \; = \; (\frac{1}{2\pi}\arg \gamma_1(x),\dots,\frac{1}{2\pi}\arg \gamma_d(x)) \; \in \;\mathbb{R}^d/\mathbb{Z}^d.\] If $v(x)$ and $v(x')$ both lie in some translate $a + S_{\delta}$ then $x - x' \in B_{\Gamma,\delta}$, and so for fixed $x' \in v(G) \cap (a + S_{\delta})$ the map $x \mapsto x - x'$ defines an injection from $v(G) \cap (a + S_{\delta})$ to $B_{\Gamma,\delta}$. Hence
\begin{equation}\label{eq476} |v(G) \cap (a + S_{\delta})| \; \leq \; |B_{\Gamma,\delta}|.\end{equation}
\textit{Proof of} (i). By a simple averaging there is some translate $a + S_{\delta}$ such that 
\[ |v(G) \cap (a + S_{\delta})| \; \geq \; |S_{\delta}||v(G)| \; = \; |S_{\delta}| \; = \; \delta^d|G|.\] The result is now immediate from \eqref{eq476}.\\[11pt]
\textit{Proof of} (ii). From \eqref{eq476} one has
\[ \left|v(B_{\Gamma,2 \delta}) \cap (a + S_{\delta})\right| \; \leq \;\left|v(G) \cap (a + S_{\delta})\right| \; \leq \; |B_{\Gamma,\delta}|.\]
Now $v(B_{\Gamma,2\delta}) \cap (a + S_{\delta})$ is empty unless $a \in S_{5\delta}$, and so
\begin{equation}\boxeq
|B_{\Gamma,2\delta}| \; = \; \frac{1}{|S_{\delta}|} \int_{\mathbb{R}^d/\mathbb{Z}^d} \left|v(B_{\Gamma,2 \delta}) \cap (a + S_{\delta})\right| da \; \leq \; |B_{\Gamma,\delta}| \cdot \frac{|S_{5\delta}|}{|S_{\delta}|} \; = \; 5^d \cdot |B_{\Gamma,\delta}|.\end{equation}

\begin{lemma}\label{lem4.2} Let $G$ be an abelian group of size $N$, let $\Gamma = \{\gamma_1,\dots,\gamma_d\}$ be a set of $d$ characters on $G$ and let $\delta > 0$. Define the smoothed Bohr neighbourhood $\SBA$ by setting $\SBA(x)  =  \int^{\infty}_0 B_{\Gamma,t}(x) \frac{e^{-t/\delta}}{\delta} \, dt$ and define $\sba(x) = \SBA(x)/\Vert \SBA \Vert_1$. Then\\
\emph{(i)} $\Vert \sba \Vert_1 = 1$;\\ 
\emph{(ii)} $\Vert \sba \Vert_{\infty} \leq 3/\delta^d N$;\\
\emph{(iii)} For all $x,y \in G$, $\left|\sba(x) - \sba(x-y)\right| \; \leq \; 5\sinh\left(\frac{\Vert y \Vert_{\Gamma}}{\delta}\right)\sba(x)$.\\ 
\emph{(iv)} For all $\eta > 0$, $\displaystyle \sum_{\Vert x \Vert_{\Gamma} \geq \eta} \sba(x) \; \leq \; 2\cdot 5^de^{-\eta/2\delta}$.
\end{lemma}
\proof (i) is trivial. Before proving (ii), note that for all $x \in G$ one has $B_{\Gamma,\delta}(x) \leq e\SBA(x)$. Indeed if $x \in B_{\Gamma,\delta}$ then $B_{\Gamma,t}(x) = 1$ for all $t \geq \delta$, and so
\[ \SBA(x) \; \geq \; \int^{\infty}_{\delta} \frac{e^{-t/\delta}}{\delta} \, dt \; = \; 1/e.\]
Part (ii) of the lemma is an immediate consequence of this observation and Lemma \ref{lem4.1} (i). To prove (iii), it is easiest to prove the corresponding statement for the unnormalised functions $\SBA$. Write $\eta = \Vert y \Vert_{\Gamma}$. Then $y \in B_{\Gamma,\eta}$.
Suppose that $B_{\Gamma,t}(x) - B_{\Gamma,t}(x - y) \neq 0$. Then either $x \in B_{\Gamma,t}$ and $x - y \notin B_{\Gamma,t}$, which means that $x \in B_{\Gamma,t} \setminus B_{\Gamma,t - \eta}$, or else $x \notin B_{\Gamma,t}$ and $x - y \in B_{\Gamma,t}$, in which case $x \in B_{\Gamma,t + \eta} \setminus B_{\Gamma,t}$. Thus certainly $x \in B_{\Gamma,t + \eta} \setminus B_{\Gamma,t - \eta}$. We have, then,
\begin{eqnarray*} \left|\SBA(x) - \SBA(x-y)\right|  & \leq &
\int^{\infty}_0 \left| B_{\Gamma,t}(x) - B_{\Gamma,t}(x - y)\right| \frac{e^{-t/\delta}}{\delta} \, dt\\ & \leq & \int^{\infty}_0 \left(B_{\Gamma,t+\eta}(x) - B_{\Gamma,t-\eta}(x)\right)\frac{e^{-t/\delta}}{\delta} \, dt \\ & = & e^{-\eta/\delta}\int^{\eta}_0 B_{\Gamma,u}(x)\frac{e^{-u/\delta}}{\delta} du + 2\sinh \left(\frac{\eta}{\delta}\right)\int^{\infty}_{\eta} B_{\Gamma,u}(x) \frac{e^{-u/\delta}}{\delta}du. \end{eqnarray*}
Now if $\eta \geq \delta$ then $e^{-\eta/\delta} \leq \sinh(\eta/\delta)$, and the result is immediate. If $\eta \leq \delta$ then we instead use the estimate
\[ e^{-\eta/\delta}\int^{\eta}_0 B_{\Gamma,u}(x)\frac{e^{-u/\delta}}{\delta} du \; \leq \; \frac{1}{\delta}\int^{\eta}_0 B_{\Gamma,u}(x)\, du \; \leq \; \frac{\eta}{\delta}B_{\Gamma,\delta}(x) \; \leq \; \frac{3\eta}{\delta}\SBA(x).\]
Since $x \leq \sinh x$ for $x \leq 1$, part (iii) of the lemma follows.
Finally we prove (iv) by using Lemma \ref{lem4.1} (ii), working once again with the unnormalised functions $\SBA$. One has
\begin{eqnarray*} \sum_{\Vert x \Vert_{\Gamma} \geq \eta} \SBA(x) & \leq & \int^{\infty}_{\eta} |B_{\Gamma,t}|\frac{e^{-t/\delta}}{\delta} \, dt \\ & = & 2\int^{\infty}_{\eta/2} |B_{\Gamma,2u}|\frac{e^{-2u/\delta}}{\delta}\,du \\ & \leq & 2\cdot 5^d \int^{\infty}_{\eta/2} |B_{\Gamma,u}| \frac{e^{-2u/\delta}}{\delta} \, du \\ & \leq & 2\cdot 5^de^{-\eta/2\delta} \int^{\infty}_{\eta/2} |B_{\Gamma,u}| \frac{e^{-u/\delta}}{\delta} \, du \\ & \leq & 2 \cdot 5^d e^{-\eta/2\delta}\Vert \SBA \Vert_1.\end{eqnarray*} This concludes the proof of (iv) and hence of Lemma \ref{lem4.2}.\endproof\\[11pt]
We are now ready for an important definition.
\begin{definition}\label{def4.5}
Let $G$ be an abelian group, let $\Gamma = \{\gamma_1,\dots,\gamma_d\}$ be a set of $d$ characters on $G$ and let $\delta > 0$. Define
\[ \psi_{\Gamma,\delta} \; = \; \sba \ast \sba,\] where $\sba$ is the normalised and smoothed Bohr neighbourhood defined in the statement of Lemma \ref{lem4.2}. 
\end{definition}
\noindent The following is a very long and rather disparate collection of properties enjoyed by the functions $\psa$, all of which will be required later on.
\begin{lemma}\label{lem4.5}
Let $\delta,\delta' > 0$ and suppose that $\Gamma,\Gamma'$ are two sets of characters with $|\Gamma| = d$, $|\Gamma| = d'$ and $\Gamma \subseteq \Gamma'$. Let $x,y$ be elements of $G$, let $f : G \rightarrow \mathbb{R}$ be a function with $\Vert f \Vert_{\infty} \leq 1$ and let $\tau \in (0,1/4)$. Consider the functions $\psa$ and $\psb$ as defined above. \\
\emph{(i)} $\psa$ has real and positive Fourier transform.\\
\emph{(ii)} $\Vert \psa \Vert_1 \; = \; 1$.\\
\emph{(iii)} $\Vert \psa \Vert_{\infty} \; \leq \; 3/\delta^d N$.\\
\emph{(iv)} Suppose that $\delta \leq 2^{-12}\tau^2/d$ and that $\gamma \in \Gamma$. Then $\Vert (\gamma - 1)\psa\Vert_1 \leq \tau$, and consequently $\widehat{\psa}(\gamma) \geq 1 - \tau$.\\
\emph{(v)} $\left|\psa(x) - \psa(x-y)\right| \leq  5\sinh\left(\frac{\Vert y \Vert_{\Gamma}}{\delta}\right)\psa(x)$.\\
For parts \emph{(vi)} -- \emph{(viii)}, assume that $\delta' \leq 2^{-13}\delta\tau^2/d'$.\\
\emph{(vi)} Let $m$ be a positive integer. Then for any $x \in G$ we have
\[ \left| (\psa^{1/2} \ast \psb \ast \stackrel{m}{\dots} \ast \psb - \psa^{1/2})(x)\right| \; \leq \; (2^m - 1)\tau \psa^{1/2}(x),\] where the notation indicates that there are $m$ copies of $\psb$ in the convolution.\\
\emph{(vii)} $\Vert \psa \ast \psb - \psa\Vert_1 \leq \tau$.\\
\emph{(viii)} $\Vert(f\psa^{1/2}) \ast \psb - \psa^{1/2} (f \ast \psb)\Vert_2 \leq \tau$.\\
\emph{(ix)} Let $\kappa,\omega > 0$, and suppose that $\gamma \in G^{\ast}$ is such that $\widehat{\psi_{\Gamma,\delta}}(\gamma) \geq \kappa$. Suppose that $\delta' \leq \omega^2 \kappa^2\delta/2^{13} d'$. Then $\Vert (\gamma - 1)\psb\Vert_1 \leq \omega$. In particular, $\widehat{\psb}(\gamma) \geq 1 - \omega$.
\end{lemma}
\proof (i) and (ii) are immediate, and (iii) is an easy consequence of (ii) and Lemma \ref{lem4.2} part (ii). To proceed further, we need to estimate the tails of $\psa$. Let $\eta > 0$ be arbitrary. We have \begin{eqnarray*} \sum_{\Vert x \Vert_{\Gamma} \geq \eta} \psa(x) & = & \sum_y \sba(y)\sum_{\Vert x \Vert_{\Gamma} \geq \eta} \sba(x - y) \\ & \leq & \sum_{\Vert y \Vert_{\Gamma} \geq \eta/2} \sba(y) + \sup_{\Vert y \Vert_{\Gamma} \leq \eta/2}\sum_{\Vert x \Vert_{\Gamma} \geq \eta} \sba(x - y).\end{eqnarray*}
Now if $\Vert x \Vert_{\Gamma} \geq \eta$ and $\Vert y \Vert_{\Gamma} \leq \eta/2$ then $\Vert x - y \Vert_{\Gamma} \geq \eta/2$. Therefore 
\begin{equation}\label{lem4.6} \sum_{\Vert x \Vert_{\Gamma} \geq \eta} \psa(x) \; \leq \; 2 \sum_{\Vert y \Vert_{\Gamma} \geq \eta/2} \sba(y)\; \leq \; 4\cdot 5^de^{-\eta/4\delta}, \end{equation} this last step following from Lemma \ref{lem4.2} (iv). Equation \eqref{lem4.6} is one that will be of much service in the sequel. To prove (iv), set $\eta = 16\delta d + 64\delta\tau^{-1}$ and note that the condition $\delta \leq 2^{-12}\tau^2/d$ implies that $\eta \leq \tau/16$. Therefore
\begin{eqnarray*}
\sum_x |1 - \gamma(x)|\psa(x) & \leq & 2\sum_{\Vert x \Vert_{\Gamma} \geq \eta} \psa(x) + \sup_{\Vert x \Vert_{\Gamma} \leq \eta} |1 - \gamma(x)| \\ & \leq & 8 \cdot 5^d e^{-\eta/4\delta} + 8\eta \\ & \leq & \tau/2 +\tau/2 \; \leq \; \tau.\end{eqnarray*}
To prove (v) write $\eta = \Vert y \Vert_{\Gamma}$, so that $y \in B_{\Gamma,\eta}$. Lemma \ref{lem4.2} (iii) tells us that $|\sba(x) - \sba(x - y)| \leq 5\sinh(\eta/\delta)\sba(x)$.  Thus
\begin{eqnarray*} |\psa(x) - \psa(x - y)| & \leq & \sum_z \sba(z)\left|\sba(x - z) - \sba(x - y - z)\right| \\ & \leq & 5\sinh(\eta/\delta)\sum_z \sba(z)\sba(x-z)\\ & = & 5\sinh(\eta/\delta)\psa(x),\end{eqnarray*} which is exactly (v).
An immediate consequence of (v) together with the inequality $|a^{1/2} - b^{1/2}| \leq a^{-1/2}|a - b|$ is the bound
\begin{equation}\label{eq100} \left|\psa^{1/2}(x) - \psa^{1/2}(x-y)\right| \; \leq \; 5\sinh\left(\frac{\Vert y \Vert_{\Gamma}}{\delta}\right)\psa^{1/2}(x),\end{equation} which will be of some use later on in the proof of the lemma. 
Now recall that parts (vi) -- (vii) of the lemma are to be proved under the assumption that $\delta' \leq 2^{-13}\delta \tau^2/d'$. We begin by estimating the sum $\sum_y 5\sinh\left(\frac{\Vert y \Vert_{\Gamma}}{\delta}\right)\psb(y)$, which arises in applications of (v) and equation \eqref{eq100} above.
Let $\eta = 160\delta'd'/\tau$, and split the sum into the ranges $\Vert y \Vert_{\Gamma} \leq \eta$ and $\Vert y \Vert_{\Gamma} \geq \eta$. The sum over the first range is trivially bounded by $\sinh(\eta/\delta)$. To bound the sum over the second range, observe that
\[ \sum_{\Vert y \Vert_{\Gamma} \geq \eta} \sinh \left(\frac{\Vert y \Vert_{\Gamma}}{\delta}\right) \psb(y) \; = \; \sinh\left(\frac{\eta}{\delta}\right) + \frac{1}{\delta} \int^{\infty}_{\eta} \cosh\left(\frac{t}{\delta}\right)\sum_{\Vert y \Vert_{\Gamma} \geq t} \psb(y) \, dt.\]
Now since $\Gamma \subseteq \Gamma'$ the set $\{y : \Vert y \Vert_{\Gamma} \geq t\}$ is a subset of $\{y : \Vert y \Vert_{\Gamma'} \geq t\}$. Moreover, since $\cosh x \leq e^x$, we can use \eqref{lem4.6} to bound 
\[ \frac{1}{\delta} \int^{\infty}_{\eta} \cosh\left(\frac{t}{\delta}\right)\sum_{\Vert y \Vert_{\Gamma'} \geq t} \psb(y) \, dt \; \leq \; \frac{4\cdot 5^{d'}}{\delta} \int^{\infty}_{\eta} e^{t/\delta - t/4\delta'}\,dt \; \leq \; \frac{16 \cdot 5^{d'} \delta'}{\delta}e^{-\eta/8\delta'}.\]
This, it can be checked, is at most $\tau/10$. It remains to observe that $2\sinh(\eta/\delta) \leq 4\eta/\delta \leq \tau/10$. Adding everything together gives the bound 
\begin{equation}\label{eq101}
\sum_y 5\sinh\left(\frac{\Vert y \Vert_{\Gamma}}{\delta}\right)\psb(y) \; \leq \; \tau.
\end{equation}
Equations \eqref{eq100} and \eqref{eq101} immediately imply that \begin{equation}\label{eq102} \sum_y \left|\psa^{1/2}(x) - \psa^{1/2}(x - y)\right|\psb(y) \; \leq \; \tau \psa^{1/2}(x),\end{equation} which easily implies part (vi) of the lemma in the case $m = 1$. To prove the result for all $m$ we proceed by induction. Supposing the result to have been proved from $m = l - 1$, we have the inequalities
\begin{eqnarray*}
& & \left| (\psa^{1/2} \ast \psb \ast \stackrel{l}{\dots} \ast \psb - \psa^{1/2})(x)\right|\\ & &  \qquad \leq \; \left| (\psa^{1/2} \ast \psb \ast \stackrel{l-1}{\dots} \ast \psb - \psa^{1/2}) \ast \psb(x)\right| + \left|(\psa^{1/2} \ast \psb - \psa^{1/2})(x)\right| \\ & & \qquad \leq \; \tau(2^{l-1} - 1)\psa^{1/2} \ast \psb(x) + \tau \psa^{1/2}(x) \\ &  & \qquad \leq \; \tau 2^{l-1} \psa^{1/2}(x) + \tau(2^{l-1} - 1)\left|(\psa^{1/2} \ast \psb - \psa^{1/2})(x)\right| \\ & & \qquad \leq \; (2^l - 1)\tau \psa^{1/2}(x).\end{eqnarray*} Thus the result is true for $m = l$ as well, which confirms (vi). Part (vii) of the lemma is an immediate consequence of (v) and equation \eqref{eq101}. Moving on to part (viii), an immediate consequence of \eqref{eq102} is
that \[ \left| (f\psa^{1/2}) \ast \psb(x) - \psa^{1/2} (f \ast \psb)(x) \right| \; \leq \; \tau\psa^{1/2}(x),\] which implies the required result. Finally, we prove statement (ix).
For any fixed $y$ we have, by (v), 
\[
\kappa |1 - \gamma(y)| \; \leq \; |\widehat{\psa}(\gamma)||1 - \overline{\gamma(y)}| \; = \; \left| \sum_x \left( \psa(x+y) - \psa(x)\right)\gamma(x)\right| \; \leq \; 5\sinh\left(\frac{\Vert y \Vert_{\Gamma}}{\delta}\right).\]
The result therefore holds if \eqref{eq101} is true with $\tau = \kappa\omega$, which it is if $\delta' \leq \kappa^2\omega^2\delta/2^{13}d'$ by a calculation almost identical to the one we did earlier to establish \eqref{eq101}.
\endproof

\section{The regularity lemma}\label{sec5} \noindent Let $R \subseteq G^{\ast}$ be a set of $d$ characters, and let $\eta \in (0,1)$. Given $R$ and $\eta$, we will always write $\psi_1 = \psi_{R,\eta_1}$, $\psi_2 = \psi_{R,\eta_2}$ where $\eta_1 = \eta$ and $\eta_2 = 2^{-40}\epsilon^6\eta/dk^4$. The functions $\psa$ are those defined in \S \ref{sec4}, whose properties were laid out in Lemma \ref{lem4.5}. Recall that in \S \ref{sec3} we outlined the need for a pair of Bohr neighbourhoods $B$ and $B'$, with $B'$ much smaller than $B$. This pair of functions $\psi_1$ and $\psi_2$ turn out to be the correct way to realise this idea. Observe that, by Lemma \ref{lem4.5} (vii), we have
\begin{equation}\label{eq10} \Vert \psi_1 \ast \psi_2 - \psi_1 \Vert_1 \; \leq \; 2^{-12}k^{-2}\epsilon^3.\end{equation}
This is the most useful way of quantifying the need that $|B + B'| \approx |B|$.
We will be working with a set called $A$, and also with sets named $A_1,\dots,A_k$.  Write $\alpha_1^x = A \ast \psi_1(x)$ and $\alpha_{i,1}^x = A_i \ast \psi_1(x)$, and similarly for $\psi_2$; it is perhaps unusual to use superscript notation for a function of $x$ like this, but it is a useful way of making some of our later formulae more readable. We will also write $A^{+x}(n) = A(x + n)$, so that $A^{+x}$ is the characteristic function of $A$ translated by $x$, and write $A_i^{+x}(n) = A_i(x + n)$.\\[11pt]
When we write $\psi_1$ and $\psi_2$ we will always assume that they come from some underlying set $R$ and parameter $\eta$ in the manner just described. We will always write $|R| = d$.
\begin{definition}[Regularity] \label{def5.1}Let $\epsilon > 0$ and let $x \in G$. We say that $x$ is an $\epsilon$-\emph{regular value} with respect to the set $A$ \emph{(}and the pair $(R,\eta)$\emph{)} if the following is true.
\begin{enumerate}
\item $\sum_{y} \left(\alpha_2^{x+y} - \alpha_1^x\right)^2\psi_1(y) \; \leq \; \epsilon^2$;
\item $\Vert \left(\left(A^{+x} - \alpha_2^x\right)\psi_2\right)^{\wedge}\Vert_{\infty} \; \leq \; \epsilon$.
\end{enumerate} We say that $(R,\eta)$ is $\epsilon$-\emph{regular} for the set $A$ if the number of $x \in G$ which are not $\epsilon$-regular is less than $\epsilon N$. 
\end{definition}
\noindent We are now in a position to state our regularity lemma. Little extra work is involved in proving a version which holds for $k$ sets $A_1,\dots,A_k$ simultaneously instead of one, and we do this in order that we may prove Theorem \ref{thm1.5}. Write $W(t)$ for a tower of twos of height $\lceil t\rceil$.
\begin{theorem}[Regularity lemma for abelian groups]\label{thm5.3} Let $k$ be a positive integer, and let $\epsilon > 0$. Let $A_1,\dots,A_k$ be subsets of $G$. Then there is a pair $(R,\eta)$ with $d \leq W(2^{11}k^2\epsilon^{-3})$ and $\eta \geq 1/W(2^{11}k^2\epsilon^{-3})$ which is $\epsilon$-regular for each $A_i$.\end{theorem}
\noindent Define the $i$th index of $(R,\eta)$, the index with respect to $A_i$, by
\begin{equation}\label{eq11} \ind_i(R,\eta) \; = \; N^{-1}\Vert \alpha_{i,1} \Vert_2^2 \; = \; N^{-1}\sum_x (\alpha_{i,1}^x)^2.\end{equation} Define also the (total) index
\[ \ind(R,\eta) \; = \; \sum_{i=1}^k \ind_i(R,\eta).\]
Observe that $\ind(R,\eta) \leq k$. The main result of this section is the following.
\begin{proposition}\label{prop5.4}
Suppose that $(R,\eta)$ is not $\epsilon$-regular for all of $A_1,\dots,A_k$. Then there is a pair $(\widetilde{R},\widetilde{\eta})$ with $|\widetilde{R}| \leq (2dk/\eta \epsilon)^{60d}$ and $\widetilde{\eta} \geq (\eta\epsilon/2dk)^{60d}$ such that $\mbox{\emph{ind}}(\widetilde{R},\widetilde{\eta}) \geq \mbox{\emph{ind}}(R,\eta) + 2^{-10}k^{-1}\epsilon^3$.
\end{proposition}
\noindent Once this is proved, it is a short step to Theorem \ref{thm5.3}. Start with the trivial pair $(R,\eta) = (\emptyset,1)$. If this is not $\epsilon$-regular for all of $A_1,\dots,A_k$ then apply Proposition \ref{prop5.4} to get a new pair $(\widetilde{R},\widetilde{\eta})$. If \textit{this} is not $\epsilon$-regular then apply Proposition \ref{prop5.4} again, and so on. The index increases by at least $2^{-10}k^{-1}\epsilon^3$ at each iteration, and so the total number of steps cannot exceed $2^{10}k^2\epsilon^{-3}$. When the algorithm finishes we have a regular pair $(R,\eta)$ and it is not hard to see that $|R|$ and $\eta$ satisfy the claimed bounds; one can afford to be incredibly crude when examining the growth of $|R|$ and the decay of $\eta$, everything other than the number of iterations being essentially irrelevant.\\[11pt]
Let us begin to address Proposition \ref{prop5.4}. Suppose that $(R,\eta)$ fails to be $\epsilon$-regular for all of $A_1,\dots,A_k$. Then there is some $i$ together with at least $\epsilon N/k$ values of $x$ which fail to be $\epsilon$-regular with respect to $A_i$. For the rest of the section write $A = A_i$; when we talk about values or pairs being regular, it will always be with respect to this underlying set $A$.\\[11pt]
Now there are two ways in which $(R,\eta)$ could fail to be $\epsilon$-regular: either clause (i) of Definition \ref{def5.1} fails for at least $\epsilon N/2k$ values of $x$, or else clause (ii) does. We shall deal with these two possibilities separately in Propositions \ref{prop5.7} and \ref{prop5.11}, which together give Proposition \ref{prop5.4} immediately. We begin with a technical lemma which will be used three times in the sequel.
\begin{lemma}\label{lem5.5} Suppose that $\phi_1,\phi_2$ and $f$ are functions from $G$ to $\mathbb{R}$ such that $\Vert \phi_1 \ast \phi_2 - \phi_1\Vert_1 = \kappa$ and $\Vert f \Vert_{\infty} \leq 1$. Write $f_i = f \ast \phi_i$. Then
\[ \Vert f_2\Vert_2^2 - \Vert f_1 \Vert_2^2 \; \geq \; \sum_x\sum_y\left(f_2(x+y) - f_1(x)\right)^2\phi_1(y) - 8\kappa N.\]
\end{lemma}
\proof It is straightforward to check the identity
\begin{equation}\label{eq832} \sum_y \left(f_2(x+y) - (f_2 \ast \phi_1)(x)\right)^2 \phi_1(y) \; = \; f_2^2 \ast \phi_1(x) - (f_2 \ast \phi_1(x))^2.\end{equation}
Write $E(x) = (f_2 \ast \phi_1 - f_1)(x)$. We have
\[ |E(x)| \; = \; \Vert f \ast \phi_2\ast \phi_1 - f \ast \phi_1\Vert_{\infty} \; \leq \; \Vert f \Vert_{\infty} \Vert \phi_2 \ast \phi_1 - \phi_1 \Vert_1 \; \leq \; \kappa.\]
Substituting $f_1 = f_2 \ast \phi_1 + E$ in two places in \eqref{eq832}, conducting some simple manipulations, and summing over $x$ proves the lemma.\endproof
\begin{corollary}\label{cor5.6} 
For any $j \in \{1,\dots,k\}$ we have $\mbox{\emph{ind}}_j(R,\eta_2) \geq \mbox{\emph{ind}}_j(R,\eta) - 2^{-9}k^{-2}\epsilon^3$.
\end{corollary}
\proof Apply lemma \ref{lem5.5} with $\phi_1 = \psi_1$, $\phi_2 = \psi_2$ and $f = A_j$. By \eqref{eq10} the hypotheses of the lemma apply with $\kappa = 2^{-12}k^{-2}\epsilon^3$. Thus
\[ \ind_j(R,\eta_2) - \ind_j(R,\eta) \; = \;  \frac{\Vert f_2\Vert_2^2 - \Vert f_1 \Vert_2^2}{N} \; \geq \;  \frac{\sum_x\sum_y\left(\alpha_2^{x+y} - \alpha_1^x\right)^2\phi_1(y)}{N} - 2^{-9}k^{-2}\epsilon^3,\] which is at least $-2^{-9}k^{-2}\epsilon^3$.\endproof
\begin{proposition}\label{prop5.7} Suppose that there are at least $\epsilon N/2k$ values of $x$ for which $\sum_y (\alpha_2^{x+y} - \alpha_1^x)^2\psi_1(y) > \epsilon^2$. Then $\mbox{\emph{ind}}(R,\eta_2) \geq \mbox{\emph{ind}}(R,\eta) + \epsilon^3/8k$.
\end{proposition}
\proof Apply Lemma \ref{lem5.5} with $\psi_1 = \phi_1$, $\psi_2 = \phi_2$ and $f = A = A_i$. The hypotheses of that lemma hold with $\kappa = \epsilon^3/32k$. One therefore has
\[ \ind_i(R,\eta_2) - \ind_i(R,\eta) \; = \;  \frac{\Vert f_2\Vert_2^2 - \Vert f_1 \Vert_2^2}{ N}  \; \geq \;  \frac{\sum_x\sum_y\left(\alpha_2^{x+y} - \alpha_1^x\right)^2\phi_1(y)}{ N} - \epsilon^3/4k \; \geq \; \epsilon^3/4k.\] It follows from this and Corollary \ref{cor5.6} that
\begin{eqnarray*} \ind(R,\eta_2) - \ind(R,\eta) & = & \sum_{j=1}^k \left(\ind_j(R,\eta_2) - \ind_j(R,\eta)\right) \\ & \geq & \frac{\epsilon^3}{4k} + \sum_{j \neq i} \left(\ind_j(R,\eta_2) - \ind_j(R,\eta)\right),\end{eqnarray*}
which is at least $\epsilon^3/8k$.\endproof\\[11pt]
We now begin working towards Proposition \ref{prop5.11}, which deals with the possibility that there are at least $\epsilon  N/2k$ values of $x$ for which $\Vert \left((A^{+x} - \alpha_2^x)\psi_2\right)^{\wedge}\Vert_{\infty} > \epsilon$. The next few lemmas are all relevant to this endeavour. Before formulating them, let us define the (approximate) orthogonal complement of $\psi_2$, $\psi_2^{\perp}$, to be the set of all $\gamma$ for which $\widehat{\psi}_2(\gamma) \geq \epsilon/6$.
\begin{lemma}\label{lem838} Suppose that $\gamma \notin \psi_2^{\perp}$, and that $\left|\left((A^{+x} - \alpha_2^x)\psi_2\right)^{\wedge}(\gamma)\right| \geq \epsilon$. Then for all $y \in B_{R,\epsilon\eta_2/60} + x$ we have $\left|\left((A^{+y} - \alpha_2^y)\psi_2\right)^{\wedge}(\gamma)\right| \geq \epsilon/2$.
\end{lemma}
\proof Since $\gamma \in \psi_2^{\perp}$, one knows that $|(A^{+x}\psi_2)^{\wedge}(\gamma)| \geq 5\epsilon/6$. Now suppose that $m \in B_{R,\epsilon^2/60}$. We have
\[ (A^{+(x+m)}\psi_2)^{\wedge}(\gamma) \; = \; \gamma(m) \sum_n A(x+n)\psi_2(n+m)\gamma(n).\]
But
\[ \sum_n A(x+n)\psi_2(n+m)\gamma(n) - (A^{+x}\psi_2)^{\wedge}(\gamma) \; = \; \sum_n A(x+n)\gamma(n)\left(\psi_2(n+m) - \psi_2(n)\right)\]
and so
\begin{eqnarray*}
\left|(A^{+(x+m)}\psi_2)^{\wedge}(\gamma) - \gamma(m)(A^{+x}\psi_2)^{\wedge}(\gamma)\right| & \leq & \sum_n \left|\psi_2(n+m) - \psi_2(n)\right| \\ & \leq & 5\sinh(\epsilon/60) \; \leq \; \epsilon/6,\end{eqnarray*} the last step being a consequence of Lemma \ref{lem4.5} (v).
Therefore $|(A^{+(x+m)}\psi_2)^{\wedge}(\gamma)| \geq 2\epsilon/3$, and so finally
\begin{equation}\boxeq \left| \left((A^{+(x+m)} - \alpha_2^{x+m})\psi_2\right)^{\wedge}(\gamma)\right| \; \geq \; \epsilon/2.\end{equation}
\begin{lemma}\label{lem5.9} Let $U \subseteq G$ be a set and let $\kappa > 0$. Then there is $K \leq (2/\kappa)^d$, disjoint sets $S_1,\dots,S_K \subseteq U$ and points $z_1,\dots,z_K \in U$, such that $\left|\bigcup_{i = 1}^K S_i \right| \geq |U|/2$ and $S_i \subseteq B_{R,\kappa} + z_i$ for each $i$.
\end{lemma}
\proof Set $\Lambda = B_{R,\kappa/2}$. We define $S_1,S_2,\dots$ and $z_1,z_2,\dots$ inductively. Suppose we have defined $S_1,\dots,S_j$, and write $U_j = U \setminus \bigcup_{i=1}^j S_i$. If $|U_j| \leq |U|/2$ then stop; at such a point one does indeed have $|\bigcup_{i=1}^j S_i| \geq |U|/2$. Otherwise, a simple averaging argument shows that there is $z$ such that $|U_j \cap (\Lambda + z)| \geq |U_j||\Lambda|/ N \geq |U_j||\Lambda|/2 N$. Set $S_{j+1} = U_j \cap (\Lambda + z)$, and let $z_{j+1}$ be any element of $S_{j+1}$. It is clear that $S_{j+1} \subseteq B_{R,\kappa} + z_{j+1}$. Now at each step of this iteration the size of $U$ is depleted by at least $|U||\Lambda|/2 N$. The maximum possible number of steps is thus no more than $ N/|\Lambda|$ which, by Lemma \ref{lem4.1} (i), is at most $(2/\kappa)^d$.\endproof
\begin{lemma}\label{lem5.10} Suppose that $\Vert\left( (A_x - \alpha_2^x)\psi_2\right)^{\wedge}\Vert_{\infty} > \epsilon$ for at least $\epsilon N/2k$ values of $x$. Then there is a pair $(\widetilde{R},\widetilde{\eta})$ with $|\widetilde{R}| \leq (2dk/\eta\epsilon)^{60d}$, $\widetilde{\eta} \geq (2dk/\eta\epsilon)^{-60d}$ and so that the associated function $\widetilde{\psi_1} = \psi_{\widetilde{R},\widetilde{\eta}}$ satisfies \\
\emph{(i)} $\Vert \psi_2 \ast \widetilde{\psi}_1 - \psi_2\Vert_1 \; \leq \; 2^{-11}k^{-2}\epsilon^3$;\\
\emph{(ii)}  There is a function $\theta : G \rightarrow G^{\ast}$ and a set $X \subseteq G$ with cardinality at least $\epsilon  N/8k$ such that for all $x \in X$ we have \[ \Vert (\theta_x - 1)\widetilde{\psi}_1\Vert_1 \; \leq \; \epsilon/8 \qquad \mbox{and} \qquad \left|\left((A^{+x} - \alpha_2^x)\psi_2\right)^{\wedge}(\theta_x)\right| \; \geq \; \epsilon/2.\]
\end{lemma}
\proof Let $Z$ be the set of all $x$ for which $\Vert \left((A_x - \alpha_2^x)\psi_2\right)^{\wedge}\Vert_{\infty} > \epsilon$ and let $\omega : Z \rightarrow G^{\ast}$ be a function such that 
\[ \left|\left((A_x - \alpha_2^x)\psi_2\right)^{\wedge}(\omega_x)\right| > \epsilon\] for all $x \in Z$. If there are at least $|Z|/2$ values of $x \in Z$ for which $\omega_x \in \psi_2^{\perp}$ then let $X$ be the set of such values, let $\widetilde{R} = R$, $\widetilde{\eta} = 2^{-80}\epsilon^{12}\eta/d^2k^8$ and let $\theta = \omega$. It is a straightforward matter to check, using Lemma \ref{lem4.5} (vii) and (ix), that both conditions (i) and (ii) are satisfied.\vspace{11pt}

\noindent Alternatively, suppose that there are at least $|Z|/2$ values of $x \in Z$ such that $\omega_x \notin \psi_2^{\perp}$. Let $U$ be the set of such points. Apply Lemma \ref{lem5.9} with $\kappa = \epsilon \eta_2/60$. This gives sets $S_1,\dots,S_K \subseteq U$ and points $z_1,\dots,z_K \in U$ where $K$, it can be checked, satisfies $K \leq (2dk/\eta\epsilon)^{50d}$. Write $\Omega = \{ \omega_{z_1},\dots,\omega_{z_K}\}$ and let $X = \bigcup S_i$. Then $|X| \geq \epsilon N/8k$. Furthermore if $x \in X$ then there is some $i$ such that $x \in B_{R,\epsilon\eta_2/60} + z_i$. By Lemma \ref{lem5.9}, this means that
\[ \left|\left((A_x - \alpha_2^x)\psi_2\right)^{\wedge}(\omega_{z_i})\right| \; \geq \;  \epsilon/2.\]
Now let $\widetilde{R} = R \cup \Omega$ and set $\widetilde{\eta} = 2^{-50}\epsilon^6\eta_2/\widetilde{d}k^4$, where $\widetilde{d} = |\widetilde{R}|$. Property (i) is a consequence of Lemma \ref{lem4.5} (vii), and (ii) follows from Lemma \ref{lem4.5} (iv). The demonstration of Lemma \ref{lem5.10} is concluded by a slightly tedious computation, which is necessary to confirm that $|\widetilde{R}|$ and $\widetilde{\eta}$ satisfy the stated bounds.\endproof
\begin{proposition}\label{prop5.11} Suppose that $\Vert\left( (A_x - \alpha_2^x)\psi_2\right)^{\wedge}\Vert_{\infty} > \epsilon$ for at least $\epsilon N/2k$ values of $x$. Then there is a pair $(\widetilde{R},\widetilde{\eta})$ with $|\widetilde{R}| \leq (2dk/\eta\epsilon)^{60d}$, $\widetilde{\eta} \geq (2dk/\eta\epsilon)^{-60d}$ and $\mbox{\emph{ind}}(\widetilde{R},\widetilde{\eta}) \geq \mbox{\emph{ind}}(R,\eta) + 2^{-10}k^{-1}\epsilon^3$.\end{proposition}
\proof Let $(\widetilde{R},\widetilde{\eta})$ be the pair constructed in Lemma \ref{lem5.10}, and let $X$ and $\theta$ be the objects associated with it so that conditions (i) and (ii) of that lemma are satisfied. Thus $|X| \geq \epsilon  N/8$ and, for all $x \in X$, one has $\left|\left((A_x - \alpha_2^x)\psi_2\right)^{\wedge}(\theta_x)\right| \geq  \epsilon/2$ and $\Vert (\theta_x - 1)\widetilde{\psi}_1\Vert_1  \leq  \epsilon/8$. We will show that if $x \in X$ then
\begin{equation}\label{eq377} \sum_y \left(\widetilde{\alpha}_1^{x+y} - \alpha_2^x\right)^2\psi_2(y) \; \geq \; \epsilon^2/16k.\end{equation} Once this is shown, an application of Lemma \ref{lem5.5} with $\phi_1 = \psi_2$, $\phi_2 = \widetilde{\psi}_1$, $f = A = A_i$ and $\kappa = 2^{-11}k^{-2}\epsilon^3$ gives
\[ \ind_i(\widetilde{R},\widetilde{\eta}) - \ind_i(R,\eta_2) \; = \;  \frac{\Vert f_2\Vert_2^2 - \Vert f_1 \Vert_2^2}{ N} \; \geq \;  \frac{\sum_x\sum_y\left(\widetilde{\alpha}_1^{x+y} - \alpha_2^x\right)^2\psi_2(y)}{ N} - 2^{-8}k^{-1}\epsilon^3,\] which is at least $2^{-8}k^{-1}\epsilon^3$.
It follows from Corollary \ref{cor5.6} that indeed $\ind_i(\widetilde{R},\widetilde{\eta}) \geq \ind_i(R,\eta) + 2^{-9}k^{-1}\epsilon^3$ and then, by another application of Corollary \ref{cor5.6}, that $\ind(\widetilde{R},\widetilde{\eta}) \geq \ind(R,\eta) + 2^{-10}k^{-1}\epsilon^3$.\\[11pt]
It remains, of course, to prove \eqref{eq377}. Suppose, throughout what follows, that $x \in X$ and write $g(n) = A^{+x}(n) - \alpha_2^x$ and $f(n) = g(n)\theta_x(n)$. Thus
\[ \left((A^{+x} - \alpha_2^x)\psi_2\right)^{\wedge}(\theta_x) \; = \; \sum_n (A^{+x}(n) - \alpha_2^x)\psi_2(n)\theta_x(n) \; = \; f\ast \psi_2(0).\]
We have, then,
\begin{equation}\label{eq16} |f \ast \widetilde{\psi}_1 \ast \psi_2(0)| \; \geq \; |f \ast \psi_2(0)| - \Vert f \ast (\psi_2 \ast \widetilde{\psi}_1 - \psi_2)\Vert_{\infty} \; \geq \; 3\epsilon/8.\end{equation}
Recalling that $\theta_x$ is a character, so that $\theta_x(n-m) = \theta_x(n)\overline{\theta_x(m)}$, one has
\begin{eqnarray}\nonumber
f \ast \widetilde{\psi}_1 \ast \psi_2(0) & = & g\theta_x \ast \widetilde{\psi}_1 \ast \psi_2(0) \\ \nonumber & = & g \ast \overline{\theta_x}\widetilde{\psi}_1 \ast \theta_x\psi_2 (0) \\ \label{eq17} & = & g \ast \widetilde{\psi}_1 \ast \theta_x \psi_2(0) + \left(g \ast \widetilde{\psi_1}(\overline{\theta} - 1)\right) \ast \theta_x\psi_2(0),\end{eqnarray}
which we may write as $E_1 + E_2$. Now $\Vert \theta_x \psi_2\Vert_1 \leq \Vert \psi_2 \Vert_1 = 1$, and so
\[ |E_2| \; \leq \; \Vert g \ast \widetilde{\psi}_1(\overline{\theta_x} - 1)\Vert_{\infty} \; \leq \; \Vert g \Vert_{\infty} \Vert (\theta_x - 1)\widetilde{\psi}_1\Vert_1 \; \leq \; \epsilon/8.\]
Together with \eqref{eq16} and \eqref{eq17}, this implies that $|E_1| \geq \epsilon/4$. But
\begin{eqnarray*} |E_1| & \leq & \sum_y \left|(A^{+x} - \alpha_2^x) \ast \widetilde{\psi}_1(y)\right|\psi_2(y)  = \sum_y \left| \widetilde{\alpha}_1^{x + y} - \alpha_2^x\right| \psi_2(y) \\ & \leq & \bigg( \sum_y (\widetilde{\alpha}_1^{x + y} - \alpha_2^x)^2\psi_2(y)\bigg)^{1/2}\bigg(\sum_y \psi_2(y)\bigg)^{1/2}  =  \bigg( \sum_y (\widetilde{\alpha}_1^{x + y} - \alpha_2^x)^2\psi_2(y)\bigg)^{1/2}.\end{eqnarray*} This confirms \eqref{eq377}, and hence the proposition.\endproof\\[11pt]
Combining Propositions \ref{prop5.7} and \ref{prop5.11} concludes the proof of Proposition \ref{prop5.4} and hence, by the comments following the statement of the proposition, of Theorem \ref{thm5.3}.\endproof

\section{The Counting Lemma}\label{sec6}
\noindent For any $k$ functions $f_1,\dots,f_k :G \rightarrow \mathbb{R}$ write
\[ T(f_1,\dots,f_k) \; = \; \sum_{x_1 + \dots + x_k = 0} f_1(x_1)\dots f_k(x_k).\]
Observe that $T$ is a multilinear operator and that
\[ T(f_1,\dots,f_k) \; = \; \int_{\gamma} \widehat{f}_1(\gamma)\dots\widehat{f}_k(\gamma)\, d\gamma,\] where the integral is taken with respect to the normalised counting measure on $G^{\ast}$ (so that it equals the sum over all $\gamma$, divided by $N$).
Throughout this section we will assume that $A_1,\dots,A_k$ are subsets of $G$ and that $\psi_1$ and $\psi_2$ come from a pair $(R,\eta)$ which is $\epsilon$-regular for the $A_i$.
The following lemma, a generalisation of a lemma in \cite{Tao}, will be used several times later on.
\begin{lemma}\label{lem6.1}
Let $f : G \rightarrow \mathbb{R}$ be a function with $\Vert f \Vert_{\infty} \leq 1$. Then
\[ \left|T\left(\psi_1^{1/2},\psi_2,\dots,\psi_2,f\psi_1^{1/2}\right) - \sum_x f(x)\psi_1(x)\right| \; \leq \; 2^k\epsilon.\]
\end{lemma}
\proof From Lemma \ref{lem4.5} (vi) we know that $|\left(\psi_1^{1/2} \ast \psi_2 \ast \stackrel{k-2}{\dots} \ast \psi_2 - \psi_1^{1/2}\right)(x)| \leq 2^k\epsilon\psi^{1/2}(x)$. Thus
\begin{eqnarray*}
& & \left|T\left(\psi_1^{1/2},\psi_2,\dots,\psi_2,f\psi_1^{1/2}\right) - \sum_x f(x)\psi_1(x)\right| \\ & & \qquad\qquad\qquad = \;\left|\sum f(x) \psi_1^{1/2}(x)\left(\psi_1^{1/2} \ast \psi_2 \ast \stackrel{k-2}{\dots} \ast \psi_2 - \psi_1^{1/2}\right)(x)\right| \\ &  &  \qquad\qquad\qquad \leq \; 2^k\epsilon \sum_x |f(x)|\psi_1(x) \; \leq \; 2^k\epsilon.\end{eqnarray*}
\begin{proposition}[Counting Lemma]\label{prop6.2}
Suppose that $A_1,\dots,A_k \subseteq G$, and that $x_1,\dots,x_k$ are $\epsilon$-regular values with $x_1 + \dots + x_k = 0$. Then
\[ \left| T\left(A_1^{+x_1}\psi_1^{1/2},A_2^{+x_2} \psi_2,\dots,A_{k-1}^{+x_{k-1}}\psi_2,A_k^{+x_k}\psi_1^{1/2}\right) - \alpha_{1,1}^{x_1}\alpha_{2,2}^{x_2}\alpha_{3,2}^{x_3}\dots\alpha_{k-1,2}^{x_{k-1}}\alpha_{k,1}^{x_k}\right| \; \leq \; 4 \cdot 2^k\epsilon.\]
\end{proposition}
\proof The multilinearity of $T$ implies that for any $2k$ functions $f_1,\dots,f_k,f'_1,\dots,f'_k$ we have
\begin{eqnarray}\nonumber T(f'_1,f'_2,\dots,f'_{k-1},f'_k) & = & T(f_1,f_2,\dots,f_{k-1},f_k) + T(f_1,f_2,\dots,f_{k-1},f'_k - f_k) \\ \nonumber & & + T(f'_1 - f_1,f_2,\dots,f_{k-1},f'_k) \\ & & + \sum_{i = 2}^{k-1} T(f'_1,\dots,f'_{i-1},f'_i - f_i,f_{i+1},\dots,f_{k-1},f'_k),\label{newstar}\end{eqnarray} which we write as $P_1 + P_2 + P_3 + \sum_{i = 2}^{k-1} Q_i$. We will apply \eqref{newstar} with $f'_1 = \alpha_{1,1}^{x_1}\psi_1^{1/2}$, $f'_i = \alpha_{i,2}^{x_i}\psi_2$ ($i = 2,\dots,k-1$), $f'_k = \alpha_{k,1}^{x_k}\psi_1^{1/2}$, $f_1 = A_1^{+x_1}\psi_1^{1/2}$, $f_i = A_i^{+x_i}\psi_2$ ($i = 2,\dots,k-1$) and $f_k = A_k^{+x_k}\psi_1^{1/2}$. Each of the terms $P_1,P_2,P_3$ and $Q_i$ will be estimated separately.
\\[11pt]
\noindent\textit{Estimation of $P_1$.} Using the multilinearity of $T$ and Lemma \ref{lem6.1} with $f = 1$, we have $\left|P_1 - \alpha_{1,1}^{x_1}\alpha_{2,2}^{x_2}\alpha_{3,2}^{x_3}\dots\alpha_{k-1,2}^{x_{k-1}}\alpha_{k,1}^{x_k}\right| \leq 2^k\epsilon$.\\[11pt]
\noindent\textit{Estimation of $P_2$.} Again we use Lemma \ref{lem6.1}, this time with $f = A_k^{+x_k} - \alpha_{k,1}^{x_k}$. We have $|P_2| \leq 2^k\epsilon$.\\[11pt]
\noindent\textit{Estimation of $P_3$.} By the Cauchy-Schwarz inequality we have
\begin{eqnarray*} |P_3| & = & \int_{\gamma} (f'_1 - f_1)^{\wedge}(\gamma) \widehat{f}_2(\gamma) \dots \widehat{f}_{k-1}(\gamma)\widehat{f'_k}(\gamma)\, d\gamma \\ & \leq & \left(\int_{\gamma} \left|(f'_1 - f_1)^{\wedge}(\gamma)\right|^2\left|\widehat{f}_2(\gamma)\right|^2\, d\gamma\right)^{1/2}\cdot \prod_{i=3}^{k-1}\Vert \widehat{f}_i\Vert_{\infty}\cdot\Vert \widehat{f}'_k\Vert_{2}.\end{eqnarray*}
Now for each $i = 3,4,\dots,k-1$ we have $\Vert \widehat{f}_i \Vert_{\infty} \leq \Vert f_i \Vert_1 \leq 1$, and it follows from Parseval's identity that $\Vert \widehat{f}'_k \Vert_2$ is also at most 1. Our attention turns, then, to the bracketed expression. Writing $F(n) = A_1^{+x_1}(n) - \alpha_{1,1}^{x_1}$ one has, using Parseval's identity,
\begin{eqnarray*} 
\int_{\gamma} \left|(f'_1 - f_1)^{\wedge}(\gamma)\right|^2|\widehat{f}_2(\gamma)|^2\, d\gamma & \leq & 
\int_{\gamma} \big|((A_1^{+x_1} - \alpha_{1,1}^x)\psi_1^{1/2})^{\wedge}(\gamma)\big|^2|\widehat{\psi_2}(\gamma)|^2\, d\gamma \\ & = & \sum_x (F\psi_1^{1/2}) \ast \psi_2(x)^2 \; = \; \Vert (F\psi_1^{1/2}) \ast \psi_2 \Vert_2^2.\end{eqnarray*}
But, using property (ii) of regularity (Definition \ref{def5.1}) and Lemma \ref{lem4.5} (viii)  we have
\begin{eqnarray*}
\Vert (F\psi_1^{1/2}) \ast \psi_2 \Vert_2 & \leq & \Vert \psi_1^{1/2} (F \ast \psi_2) \Vert_2 + \Vert (F\psi_1^{1/2})\ast \psi_2 - \psi_1^{1/2} (F\ast\psi_2)\Vert_2 \\ & = & \big(\sum_y \psi_1(y)(\alpha_{1,1}^{x_1+y} - \alpha_{1,2}^y)^2\big)^{1/2} + \Vert (F\psi_1^{1/2})\ast \psi_2 - \psi_1^{1/2} (F\ast\psi_2)\Vert_2 \\ & \leq & 2\epsilon.\end{eqnarray*}
It follows that $|P_3| \leq 2\epsilon$.\\[11pt]
\noindent\textit{Estimation of the $Q_i$.} For each $i = 2,\dots,k-1$ the quantity $Q_i$ succumbs to the estimate
\begin{eqnarray*}
|Q_i| & = & \big|\sum_{\gamma}\widehat{f}'_1(\gamma)\dots \widehat{f}'_{i-1}(\gamma)(f'_i - f_i)^{\wedge}(\gamma)\widehat{f}_{i+1}(\gamma)\dots\widehat{f}_{k-1}(\gamma)\widehat{f}'_k(\gamma)\big| \\ & \leq & \Vert \widehat{f}'_1\Vert_2 \cdot  \prod_{j = 2}^{i-1} \Vert \widehat{f}'_j \Vert_{\infty} \cdot \left\Vert (f'_i - f_i)^{\wedge}\right\Vert_{\infty} \cdot \prod_{j = i+1}^{k-1} \Vert \widehat{f}_j\Vert_{\infty} \cdot \Vert \widehat{f}'_k \Vert_2\\ & \leq & \big\Vert (f'_i - f_i)^{\wedge}\big\Vert_{\infty}\\ & = & \left\Vert    \left((A_i^{+x_i} - \alpha_{i,2}^{x_i})\psi_2\right)^{\wedge} \right\Vert_{\infty}.\end{eqnarray*}
By property (ii) of regularity (Definition \ref{def5.1}), this is at most $\epsilon$.\\[11pt]
Combining the estimates for $P_1,P_2,P_3$ and $Q_i$ with \eqref{newstar} completes the proof of Proposition \ref{prop6.2}. \endproof

\section{$k$-tuples of sets with few zero-sums}\label{sec7} \noindent We retain the notation of the last two sections.
We are now in a position to prove Theorem \ref{thm1.5}. Theorem \ref{thm1.4} and Corollary \ref{cor1.6} will be easy consequences. Let us restate the result for the reader's convenience. If $\prod_{i=1}^k X_i$ is a cartesian product of sets in $G$, we define a \emph{zero-sum $k$-tuple} to be a $k$-tuple $(x_1,\dots,x_k) \in \prod_{i=1}^k X_i$ with $x_1 + \dots + x_k = 0$.\\[11pt]
\noindent\textbf{Theorem \ref{thm1.5}} 
\emph{Let} $k \geq 3$ \emph{be a fixed integer, and suppose that} $A_1,\dots,A_k$ \emph{are subsets of} $G$ \emph{such that} $\prod_{i=1}^k A_i$ \emph{has} $o(N^{k-1})$ \emph{zero-sum} $k$-\emph{tuples}. \emph{Then we may remove} $o(N)$ \emph{elements from each} $A_i$ \emph{so as to leave sets} $A'_i$, \emph{such that} $\prod_{i=1}^k A'_i$ \emph{has no zero-sum} $k$-\emph{tuples}.
\\[11pt]
\proof We begin by setting up a definition analogous to Definition \ref{def1.7}.
\begin{definition}\label{def7.1} Suppose that $A_1,\dots,A_k \subseteq G$ and let $(R,\eta)$ be $\epsilon$-regular for the $A_i$. We define the \emph{reduced sets} $A'_i$ to be the sets obtained by deleting, for each $i = 1,\dots,k$, all $x \in A_i$ for which $x$ is not a regular value, or for which $\alpha_{i,1}^x \leq 4\epsilon^{1/k}$ or $\alpha_{i,2}^x \leq 4\epsilon^{1/k}$.
\end{definition}
\noindent It was a simple matter to show that obtaining the reduced graph from the original graph involved the deletion of rather few edges (cf. \eqref{eq1}). To show that $|A \setminus A'|$ is small is a little subtle. In fact it is to obtain such a result that we have been dealing with the functions $\psa$ rather than the functions $\sba$. Up until now, either would have worked.
\begin{lemma}\label{lem7.2} Let $A \subseteq G$ and let $\psi = \psi_{\Gamma,\delta}$ for some choice of $\Gamma \subseteq G^{\ast}$ and some $\delta > 0$. Let $\rho > 0$. Then the number of $x \in A$ for which $A \ast \psi(x) \leq \rho$ is no more than $\rho N$.\end{lemma}
\proof Write $\beta = \beta_{\Gamma,\delta}$, so that $\psi = \beta \ast \beta$. Let $S$ be the set of all $x \in A$ for which $A \ast \psi(x) \leq \rho$. Certainly, then, $S \ast \psi (x) \leq \rho$ for all $x \in S$. Thus we have
\begin{eqnarray*}
\frac{|S|^2}{N} \; = \; \frac{1}{N}\big(\sum_x S(x)\beta(x)\big)^2 & \leq & \sum_x S \ast \beta(x)^2 \; = \; \sum_x S(x) S \ast \beta \ast \beta(x) \\ & = & \sum_x S(x) S \ast \psi (x) \; \leq \; \rho |S|.\end{eqnarray*}
The result follows immediately.\endproof\\[11pt]
Recall Definition \ref{def7.1}. Since the number of $x$ which fail to be regular is no more than $\epsilon N$, it follows that
\[ |A'_i| \; \geq \; |A_i| - 10k\epsilon^{1/k}N\] for each $i$.\endproof\\[11pt]
\noindent\textit{Proof of Theorem \ref{thm1.5}.} Suppose that there are at most $\delta N^{k-1}$ zero-summing $k$-tuples $(a_1,\dots,a_k) \in \prod_{i=1}^kA_i$. Choose a function $\epsilon = \epsilon(\delta)$ such that $\epsilon \rightarrow 0$ as $\delta \rightarrow 0$, but such that $\delta^{-1} > W(2^{12}k^2\epsilon^{-3})$. 
This means that there is a pair $(R,\eta)$ which is $\epsilon$-regular for $A$, and for which the associated constants $d = |R|$ and $\eta_2 = 2^{-40}\epsilon^6\eta/dk^4$ satisfy the condition $3^k\delta/\eta_2^{dk} < \epsilon$. Again, this is an easy check since everything but the height of the power of twos is essentially irrelevant. Consider the reduced sets $A'_i$ coming from such a regular partition $(R,\eta)$. As we have seen, $|A'_i| \geq |A_i| - 10k\epsilon^{1/k}N$. We claim that there are no zero-sum $k$-tuples in $\prod_{i=1}^k A'_i$. Indeed, suppose that $x_i \in A'_i$ are such that $x_1 + \dots + x_k = 0$. Then the counting lemma tells us that 
\[ T(A_1^{+x_1}\psi_1^{1/2},A_2^{+x_2}\psi_2,\dots,A_{k-1}^{+x_{k-1}}\psi_2,A_k^{+x_k}\psi_1^{1/2}) \; \geq \; \epsilon.\]
However every $k$-tuple $(u_1,\dots,u_k)$ contributing to the sum
\begin{eqnarray*} & &  T(A_1^{+x_1}\psi_1^{1/2},A_2^{+x_2}\psi_2,\dots,A_{k-1}^{+x_{k-1}}\psi_2,A_k^{+x_k}\psi_1^{1/2}) \\ & & \qquad \qquad = \; \sum_{u_1 + \dots + u_k = 0}A_1^{+x_1}\psi_1^{1/2}(u_1)A_2^{+x_2}\psi_2(u_2)\dots A_{k-1}^{+x_{k-1}}\psi_2(u_{k-1})A_k^{+x_k}\psi_1^{1/2}(u_k)\end{eqnarray*} leads to a zero-sum $k$-tuple $(x_1 + u_1,x_2+u_2,\dots,x_k + u_k) \in \prod_{i=1}^k A_i$. By Lemma \ref{lem4.5} (iii) both $\Vert \psi_1\Vert_{\infty}$ and $\Vert \psi_2 \Vert_{\infty}$ are at most $3/\eta_2^d N$. Since there are no more than $\delta N^{k-1}$ zero-sum $k$-tuples in $\prod_{i=1}^k A_i$, this means that
\[ T(A_1^{+x_1}\psi_1^{1/2},A_2^{+x_2}\psi_2,\dots,A_{k-1}^{+x_{k-1}}\psi_2,A_k^{+x_k}\psi_1^{1/2}) \; \leq \; \frac{3^{k}\delta}{\eta_2^{dk}} \; < \; \epsilon.\] This is a contradiction, and so $\prod_{i=1}^k A'_i$ does indeed cannot any zero-sum $k$-tuples.\endproof\\[11pt]
\noindent\textit{Proof of Theorem \ref{thm1.4}.} Simply apply Theorem \ref{thm1.5} with $k = 3$ and $A_1 = A_2 = A_3$.\endproof\\[11pt]
\noindent\textit{Proof of Corollary \ref{cor1.6}.} Apply Theorem \ref{thm1.5} with $G = \mathbb{Z}/2N\mathbb{Z}$, $k = 3$. If $A \subseteq [N]$ then $A$ may be identified with a subset of $\mathbb{Z}/2N\mathbb{Z}$ by reducing modulo $2N$. Set $A_1 = A_2 = A$ and $A_3 = - A$. The result follows by observing that if $x,y,z \in A$ then $x + y \equiv z \pmod{2N}$ if, and only if, $x + y = z$.\endproof\\[11pt]
\noindent\textit{Remarks on bounds.} Consider Corollary \ref{cor1.6} as formulated in the abstract of the paper. That is, if $A \subseteq [N]$ has $\delta N^2$ summing triples (triples with $x + y = z$) then it may be made sum-free by removing $\delta' N$ elements. Our proof gives an awful dependence between $\delta$ and $\delta'$ of the form $1/\delta = W(\delta^{\prime -C})$. One may conjecture that a much stronger result should be true. It is, however, too optimistic to hope that the dependence between $\delta'$ and $\delta$ might be polynomial, and we close this section by giving a very brief sketch of why this is so. Let $\delta > 0$, let $N$ be a large positive integer and let $p$ and $q$ be distinct primes to be chosen later. Recall that by a construction of Behrend \cite{Beh} there is a set $B \subseteq \mathbb{Z}/p\mathbb{Z}$ with $|B| \sim \exp(-C_1(\log p)^{1/2})p$ with the \textit{Behrend property}, that is the only triples $(x,y,z)$ with $x + y = 2z$ are the trivial ones for which $x = y = z$. This construction involves choosing suitable integers $r$ and $d$, taking the lattice points on the sphere $\Sigma((0,\dots,0),r) \subseteq \mathbb{R}^d$ and then using an affine transformation to project to $\mathbb{Z}/p\mathbb{Z}$. Now the lattice points on a sphere clearly have the Behrend property, but the same is also true of the union $X \cup 2X \cup 4X$, where $X$ is the set of lattice points on the sphere $\Sigma((r,0,\dots,0),r')$, provided that $r' < r/7$. Using this set instead, we may find $S \subseteq \mathbb{Z}/p\mathbb{Z}$ with $|S| \sim \exp(-C_2(\log p)^{1/2})p$ such that $S$ has both the Behrend property and also the additional property that $S \cap 2S \cap 4S$ has cardinality at least $\exp(-C_3(\log p)^{1/2})p$. Set $T = S \cap 2S$: then $|T \cap 2T| = |S \cap 2S \cap 4S|$ satisfies the lower bound just mentioned. By choosing $p \sim \exp(-C_4(\log(\frac{1}{\delta}))^{1/2})/\delta$ we may take $|T| \sim \delta p^2$. Now choose $q$ so that $pq \sim N$ and consider the set $U = T \times \mathbb{Z}/q\mathbb{Z}$ (which may be considered as a subset of $\mathbb{Z}/N\mathbb{Z}$). If $(t_1,a_1) + (t_2,a_2) = (t_3,a_3)$ then certainly $t_1 + t_2 = t_3$. But $t_3 = 2s$ for some $s \in S$ and so $t_1 + t_2 = 2s$, which implies that $t_1 = t_2$. It follows that the number of summing triples in $U$ is bounded above by $|T|q^2 \sim \delta N$. However, in order to remove all summing triples from $U$ we must delete at least one of $(x,a)$ and $(2x,2a)$ for all $x \in T \cap 2T$, $a \in \mathbb{Z}/q\mathbb{Z}$, and to do this requires the removal of at least $\exp(-C_5(\log(\frac{1}{\delta}))^{1/2})N$ elements.\vspace{11pt}

\noindent I do not know a similar example in $(\mathbb{Z}/2\mathbb{Z})^n$ -- that is, it may be that the dependence between $\delta'$ and $\delta'$ in Theorem \ref{thm1.4} is polynomial when $G = (\mathbb{Z}/2\mathbb{Z})^n$. It is my belief that the only bounds known in the classical Proposition \ref{prop1.3} are also of tower type.

\section{A question of Bergelson, Host and Kra} \label{sec7.5} \noindent We again use the notation of \S \ref{sec5} and \ref{sec6}. In this section we prove Theorems \ref{thm1.8} and \ref{thm1.9}, starting with the former. We restate it for the reader's convenience.\vspace{11pt}

\noindent\textbf{Theorem \ref{thm1.8}} 
\emph{Suppose that} $\alpha,\epsilon > 0$. \emph{Then there is} $N_0(\alpha,\epsilon)$ \emph{such that if} $G$ \emph{is an abelian group of size} $N > N_0(\alpha,\epsilon)$ \emph{with} $N$ \emph{odd, and if} $A \subseteq G$ \emph{has size} $\alpha N$, \emph{then there is some} $d \neq 0$ \emph{such that} $A$ \emph{has at least} $(\alpha^3 - \epsilon)N$ \emph{three-term arithmetic progressions with common difference} $d$.\vspace{11pt}

\proof Define three sets $A_1,A_2,A_3$ by $A_1 := A$, $A_2 := -2A = \{-2x : x \in A\}$ and $A_3 := A$. By Theorem \ref{thm5.3} we may find a pair $(R,\eta)$ which is $\epsilon$-regular for $A_1,A_2$ and $A_3$. For a given $x \in G$, set $x_1 = x_3 := x$ and $x_2 := -2x$. If $x_1,x_2$ and $x_3$ are all regular values then the counting lemma (Proposition \ref{prop6.2}) tells us that
\[ \left| T(A_1^{+x_1}\psi_1^{1/2}, A_2^{+x_2} \psi_2, A_3^{+x_3} \psi_1^{1/2}) - \alpha^{x_1}_{1,1} \alpha^{x_2}_{2,2} \alpha^{x_3}_{3,1}\right| \leq 32\epsilon.\]
Since the number of $x$ for which some $x_i$ fails to be $\epsilon$-regular for $A_i$ is at most $2\epsilon N$, it follows that
\begin{equation}\label{eq8.101} \sum_x T(A_1^{+x_1}\psi_1^{1/2}, A_2^{+x_2} \psi_2, A_3^{+x_3} \psi_1^{1/2}) \geq \sum_x \alpha^{x_1}_{1,1} \alpha^{x_2}_{2,2} \alpha^{x_3}_{3,1} - 34 \epsilon N.\end{equation}
Now it is a fairly straightforward matter to check that 
\[ \sum_x T(A_1^{+x_1}\psi_1^{1/2}, A_2^{+x_2} \psi_2, A_3^{+x_3} \psi_1^{1/2}) = \sum_d P(A;d)\nu(d),\]
where $P(A;d)$ is the number of 3-term arithmetic progressions in $A$ with common difference $d$ and the weight $\nu(d)$ is given by
\[ \nu(d) := \sum_y \psi_1^{1/2}(y) (\textstyle\frac{1}{2}\displaystyle\psi_2)(y + d) \psi_1^{1/2}(y + 2d).\] Here we have written $\frac{1}{2}\psi_2(t) = \psi_2(2t)$. Together with \eqref{eq8.101} this implies that 
\begin{equation}\label{eq8.102} \sum_d P(A;d) \nu(d) \geq \sum_x \alpha^{x_1}_{1,1} \alpha^{x_2}_{2,2} \alpha^{x_3}_{3,1} - 34 \epsilon N.\end{equation}
Now we have $\alpha_{1,1}^{x_1} = \alpha_{3,1}^{x_3} = (A \ast \psi_1)(x)$, and furthermore (since $N = |G|$ is odd) $\alpha_{2,2}^{x_2} = (A \ast \frac{1}{2}\psi_2)(x)$. For notational convenience write $\alpha_1(x) := \alpha_{1,1}^{x_1} = \alpha_{3,1}^{x_3}$ and $\alpha_2(x) := \alpha_{2,2}^{x_2}$. Now we have $\psi_1 = \beta_1 \ast \beta_1$ and $\frac{1}{2}\psi_2 = \frac{1}{2}\beta_2 \ast \frac{1}{2}\beta_2$, where $\beta_1,\beta_2$ are the smoothed Bohr cutoffs used to define $\psi_1$ and $\psi_2$. It follows that 
\begin{eqnarray*} \sum_x \alpha_1(x)\alpha_2(x) & = & \sum_{x} (A \ast \beta_1 \ast \beta_1)(x)(A \ast \textstyle\frac{1}{2}\displaystyle \beta_2 \ast \textstyle\frac{1}{2}\displaystyle \beta_2)(x) \\ & = & \sum_x (A \ast \beta_1 \ast \textstyle\frac{1}{2}\displaystyle \beta_2)(x)^2 \\ & \geq & N^{-1} (\sum_x (A \ast \beta_1 \ast \textstyle\frac{1}{2}\displaystyle \beta_2)(x) )^2  \; = \; \alpha^2 N.\end{eqnarray*}
Combining this with the fact that $\sum_x \alpha_2(x) = \alpha N$, we obtain the inequality
\[ \sum_x \alpha^{x_1}_{1,1} \alpha^{x_2}_{2,2} \alpha^{x_3}_{3,1} = \sum_x \alpha_1^2(x)\alpha_2(x) \geq \big(\sum_x \alpha_1(x)\alpha_2(x)\big)^2/\sum_x \alpha_2(x) \geq \alpha^3 N.\]
Comparing with \eqref{eq8.102} yields
\[  \sum_d P(A;d) \nu(d) \geq (\alpha^3 - 34 \epsilon) N,  \] and so provided $N > N_0(\alpha,\epsilon)$ we have
\begin{equation}\label{eq8.103} \sum_{d \neq 0} P(A;d) \nu(d) \geq (\alpha^3 - 35 \epsilon) N.\end{equation}
Note that if $N$ is too small then the cutoffs $\psi_1,\psi_2$ will be almost entirely supported at zero, and such a conclusion would not be correct.
Now it is easy to see that $\sum_d \nu(d) = T(\psi_1^{1/2},\psi_2,\psi_1^{1/2})$,
and so Lemma \ref{lem6.1} with $f = 1$ gives
\begin{equation}\label{eq8.110} \sum_d \nu(d) \leq  1 + 8\epsilon.\end{equation}
Thus there is some $d \neq 0$ such that $P(A;d) \geq (\alpha^3 - 35\epsilon)/(1 + 8\epsilon)$, which implies Theorem \ref{thm1.8} after redefining $\epsilon$.\endproof\vspace{11pt}

\noindent We move on now to outline the proof of Theorem \ref{thm1.9}, which is the same result but with $G$ replaced by $\{1,\dots,N\}$. \vspace{11pt}

\noindent\textit{Proof of Theorem \ref{thm1.9}.} Suppose that $A \subseteq \{1,\dots,N\}$ has density $\alpha$. We may regard $A$ as a set $\overline{A} \subseteq \mathbb{Z}/N\mathbb{Z}$ in a natural way; observe, however, that 3-term arithmetic progressions in $\overline{A}$ need not be three-term progressions in $A$. To get around this problem, a trick is required. As before we take a pair $(R,\eta)$ which is $\epsilon$-regular for $\overline{A}$, but now we insist that $\gamma_{1/2},\gamma_1 \in R$, where $\gamma_1$ is the character $x \mapsto e^{2\pi i x/N}$ and $\gamma_{1/2}$ maps $x$ to $e^{2\pi i \overline{2}x/N}$, where $\overline{2}$ is the multiplicative inverse of $2 \md{N}$. It is a trivial matter to achieve this (at the expense of an inconsequentially worse upper bound for $|R|$) by starting the iteration used to prove Theorem \ref{thm5.3} with the pair $(R,\eta) = (\{\gamma_1,\gamma_{1/2}\},1)$ instead of $(\emptyset,1)$.\vspace{11pt}

\noindent Now we argue exactly as in the proof of Theorem \ref{thm1.8}, obtaining the bound \eqref{eq8.103}, that is to say
\begin{equation}\label{eq8.104} \sum_d P(\overline{A};d) \nu(d) \geq (\alpha^3 - 35\epsilon)N.\end{equation}
Given $d \in \mathbb{Z}/N\mathbb{Z}$, write $|d|$ for the magnitude of that residue $\overline{d} \in -\{(N-1)/2,\dots,(N-1)/2\}$ with $\overline{d} \equiv d \md{N}$. We will show that almost all of the sum on the left in \eqref{eq8.104} is concentrated on those $d$ for which $|d|$ is small, which is good as such $d$ are rather likely to correspond to 3-term progressions in $A$, rather than just $\md{N}$ progressions in $\overline{A}$. \vspace{11pt}

\noindent Suppose then that $|d| \geq \epsilon N$. Write $S_{\delta}$ for the set all all $x$ such that $|x| \geq \delta N$. One has
\begin{eqnarray}\nonumber
\sum_{|d| \geq \epsilon N} \nu(d) & \leq & \sum_{x,d} \psi_1^{1/2}(x)(\textstyle \textstyle\frac{1}{2}\displaystyle  \displaystyle\psi_2)(x+d)\psi_1^{1/2}(x+2d)S_{\epsilon}(d) \\ \nonumber & \leq & \big(\sum_{x,d} \psi_1(x)(\textstyle \textstyle\frac{1}{2}\displaystyle  \displaystyle\psi_2)(x+d) S_{\epsilon}(d)\big)^{1/2} \big(\sum_{x,d} \psi_1(x+d)(\textstyle \textstyle\frac{1}{2}\displaystyle  \displaystyle\psi_2)(x+2d) \big)^{1/2} \\ \label{eq8.106} & = & \big(\sum_{x,d} \psi_1(x)(\textstyle \textstyle\frac{1}{2}\displaystyle  \displaystyle\psi_2)(x+d) S_{\epsilon}(d)\big)^{1/2}.\end{eqnarray}
Now we have $S_{\epsilon}(d) \leq S_{\epsilon/2}(x) + S_{\epsilon/2}(x+d)$. We split the sum in \eqref{eq8.106} into two parts accordingly, that is to say as
\[ \Sigma_1 := \sum_{x,d} \psi_1(x)(\textstyle \textstyle\frac{1}{2}\displaystyle  \displaystyle\psi_2)(x+d)S_{\epsilon/2}(x) = \sum_x \psi_1(x)S_{\epsilon/2}(x)\]
and 
\[ \Sigma_2 := \sum_{x,d} \psi_1(x)(\textstyle \textstyle\frac{1}{2}\displaystyle  \displaystyle\psi_2)(x+d)S_{\epsilon/2}(x+d) = \sum_y \psi_2(2y)S_{\epsilon/2}(y).\]
To estimate $\Sigma_1$, observe that if $|x| \geq \epsilon N/2$ then $\Vert x \Vert_R \geq \epsilon/2$, the notation being that of \S \ref{sec3}, by virtue of the fact that $\gamma_1 \in R$. If, as we may, we assume that $|R|$ is much smaller than $1/\eta$ and that $\eta$ is enormously smaller than $\epsilon$ it follows from \eqref{lem4.6} that $\Sigma_1 \leq \epsilon$ (this is true by a huge margin).\vspace{11pt}

\noindent Turning to $\Sigma_2$, note that if $|y| \geq \epsilon N/2$ then $\Vert 2y \Vert_R \geq \epsilon/2$, this following from that fact that $\gamma_{1/2} \in R$. Once again, then, we have the estimate $\Sigma_2 \leq \epsilon$ by a vast margin.\vspace{11pt}

\noindent Collating these observations together with \eqref{eq8.106} leads to the bound 
\[ \sum_{|d| \geq \epsilon N} \nu(d) \leq 2\epsilon,\] and so in view of \eqref{eq8.104} we obtain
\[ \sum_{\substack{|d| \leq \epsilon \\ d \neq 0}} P(\overline{A},d)\nu(d) \geq (\alpha^3 - 37\epsilon)N.\] Together with \eqref{eq8.110} this implies that there is some $d \neq 0$, $|d| \leq \epsilon N$, such that $P(\overline{A},d) \geq (\alpha^3 - 45\epsilon)N$. Now of the arithmetic progressions in $\overline{A}$ with common difference $d$, at most $2\epsilon N$ do not actually correspond to genuine progressions of integers under the inverse of the projection map $\{1,\dots,N\} \rightarrow \mathbb{Z}/N\mathbb{Z}$. In conclusion, then, the set $A$ contains at least $(\alpha^3 - 47 \epsilon)N$ three-term progressions with common difference $d$. \endproof

\section{Miscellaneous remarks} \label{sec8} \noindent In this section we assemble a variety of remarks concerning Theorem \ref{thm5.3}, its application, and its relationship with results in the literature.\\[11pt]
\noindent\textit{I. Relationship with Szemer\'edi's regularity lemma.} In addition to the analogies we have already drawn between Theorem \ref{thm5.3} and SzRL there is another, more formal, link between the two theorems. When applying SzRL in number theory one might consider a graph derived from a subset of an abelian group by something akin to the Cayley graph construction (cf. \cite{FGR,RuzSzem}). Perhaps the simplest situation is the following. Let $G = (\mathbb{Z}/2\mathbb{Z})^n$, let $A \subseteq G$ and let $\Gamma$ be a bipartite graph on vertex set $G \times \{0,1\}$, $(g,0)$ being joined to $(g',1)$ precisely if $g + g' \in A$. If $X \subseteq G$ we will write $X_i$ for $X \times \{i\}$ ($i = 0,1$).\vspace{11pt}

\noindent For the rest of this discussion we revert to the language of \S \ref{sec2}.
Suppose that $H \leq G$ is $\epsilon^2$-regular for $A$, this having a fairly simple meaning since $G = (\mathbb{Z}/2\mathbb{Z})^n$. We claim that if $x_1,x_2 \in G$, and if $x = x_1 - x_2$ is an $\epsilon^2$-regular value, then the pair $(H + x_1,H + x_2)$ is $\epsilon$-regular in the \textit{graph-theoretic} sense of \S \ref{sec1}. To see this, suppose that $U + x_1 \subseteq H + x_1$ and $V + x_2 \subseteq H + x_2$ both have cardinality at least $\epsilon|H|$. Then $e(U + x_1,V + x_2)$ is exactly $\sum_{u,v} A^{+x}(u + v)U(u)V(v)$. This may be estimated by Fourier techniques on $H$. Indeed
\begin{eqnarray*}
e(U + x_1,V + x_2) & = & |H|^{-1}\sum_{\gamma} \widehat{A^{+x}_H}(\gamma)\widehat{U}(\gamma)\widehat{V}(\gamma) \\ & = & \frac{|A||U||V|}{|H|} + \frac{1}{|H|}\sum_{\gamma \neq 0} \widehat{A^{+x}_H}(\gamma)\widehat{U}(\gamma)\widehat{V}(\gamma).\end{eqnarray*}
It follows that
\begin{eqnarray*}
\left| d(U + x_1,V + x_2) - d(H + x_1,H + x_2)\right| & = & \left|\frac{1}{|H||U||V|}\sum_{\gamma \neq 0}\widehat{A^{+x}_H}(\gamma)\widehat{U}(\gamma)\widehat{V}(\gamma)\right| \\ & \leq & \frac{1}{|H||U||V|}\sup_{\gamma \neq 0} \left|\widehat{A^{+x}_H}(\gamma)\right|\Vert \widehat{U} \Vert_2 \Vert \widehat{V} \Vert_2 
\\ & \leq & \frac{\epsilon^2|H|}{|U|^{1/2}|V|^{1/2}} \;\; \leq \; \epsilon.\end{eqnarray*} 
This confirms the claim. Now partition both vertex classes of $\Gamma$ into cosets $H + x_1,\dots,H + x_k$. For fixed $i$ there are at most $\epsilon^2 k$ values of $j$ for which $x_i - x_j$ is not regular, and so this partition is $\epsilon$-regular in the graph-theoretic sense of Szemer\'edi.\vspace{11pt}

\noindent Have we, then, simply recovered SzRL? In fact, rather more has been achieved. The graph $\Gamma$ was of a special type (essentially a Cayley graph over $G$) but in return we were able to insist that the vertex classes in SzRL were not arbitrary sets, but subgroups of $G$. It is hard to formulate this principle at all precisely for groups other than $(\mathbb{Z}/2\mathbb{Z})^n$. Morally speaking, however, Theorem \ref{thm5.3} says that if we have a Cayley-type graph over a group $G$ then the classes in SzRL may be chosen to have a rather strong structure which is related to $G$.\\[11pt]
\noindent\textit{II. Enumeration of sum-free and related sets.} Let $b = (b_1,\dots,b_k)$ be a fixed $k$-tuple of non-zero integers. For any real numbers $x_1,\dots,x_k$ we write $L_b(x_1,\dots,x_k) = b_1x_1 + \dots + b_kx_k$. We say that a set $A \subseteq [N]$ is strongly $L_b$-free if there are no solutions to $L_b(a_1,\dots,a_k) = 0$ with $a_i \in A$ for all $i$.
$A$ is deemed to be weakly $L_b$-free if the only solutions to $L_b(a_1,\dots,a_k) = 0$ are \textit{trivial}, that is to say they arise by partitioning $[k]$ into $I_1 \cup \dots \cup I_t$ such that $\sum_{i \in I_j} b_i = 0$ for each $j$, and then taking all of the $a_i$ $(i \in I_j)$ to be equal. When $b = (1,1,-2)$, a strongly $L$-free set must be empty, and a weakly $L$-free set is the same thing as a Behrend set (cf. \S \ref{sec7}). When $b = (1,1,-1)$, the notions of strongly and weakly $L$-free coincide with that of a sum-free set. The reader is referred to \cite{Ruz1} for more information on solving linear equations in sets of integers. \\[11pt]
In this subsection we use Theorem \ref{thm5.3} to get estimates on the number of weakly $L_b$-free subsets of $[N]$. Many of our results would extend to arbitrary abelian groups, but the discussion of general linear forms is complicated by the possibility of torsion and we do not give it here.\\[11pt]
Machinery for counting sum-free sets was developed by I.Z. Ruzsa and the author in a series of papers \cite{Gre2,GreRuz1,GreRuz2}. The following result may be proved by extending the methods used in those papers in a straightforward manner.
\begin{proposition}[Granularization]\label{prop8.1} Let $L = L_b$ be a fixed linear form as above. Then there is a family $\mathcal{F}$ of subsets of $[N]$ with the following properties:
\begin{enumerate}
\item $|\mathcal{F}| = 2^{o(N)}$;
\item If $A \subseteq [N]$ is weakly $L$-free, then $A \subseteq F$ for some $F \in \mathcal{F}$;
\item Each $F \in \mathcal{F}$ has $o(N^{k-1})$ solutions to $L(f_1,\dots,f_k) = 0$.\end{enumerate}
\end{proposition}
\noindent Now Theorem \ref{thm5.3} gives structural information about the sets in $\mathcal{F}$. Indeed suppose $F \in \mathcal{F}$ and set $M = 2(|b_1| + \dots + |b_k|)N$. Apply Theorem \ref{thm1.4} with $G = \mathbb{Z}/M\mathbb{Z}$ and $A_i = F_i\pmod{M}$, where $F_i = b_iF$. The choice of $M$ guarantees that $f_1 + \dots + f_k \pmod{M}$ if, and only if, $f_1 + \dots + f_k = 0$. It follows that there is a set $F' \subseteq F$ with $|F \setminus F'| = o(N)$ and such that $F'$ is strongly $L$-free. This leads to the following improvement of Proposition \ref{prop8.1}:
\begin{proposition}\label{prop8.2}
Let $L = L_b$ be a fixed linear form as above. Then there is a family $\mathcal{F}'$ of subsets of $[N]$ with the following properties:\begin{enumerate}
\item $|\mathcal{F}'| = 2^{o(N)}$;
\item If $A \subseteq [N]$ is weakly $L$-free, then $|A \setminus F'| = o(N)$ for some $F' \in \mathcal{F}'$;
\item Each $F' \in \mathcal{F}'$ is strongly $L$-free.\end{enumerate}
\end{proposition}
\noindent Note that this implies that $r_3(N) = o(N)$, since when $b = (1,1,-2)$ the collection $\mathcal{F}$ must consist of just the empty set.
Using Proposition \ref{prop8.2} and the methods of our papers with Ruzsa, one can prove the following. Write $\mbox{LF}(N)$ for the collection of all weakly $L$-free subsets of $[N]$. 
\begin{theorem}\label{thm8.3}
Let $L = L_b$ be a fixed linear form. For each $N$ denote by $f_L(N)$ be the cardinality of the largest strongly $L$-free subset of $[N]$. Then $\log_2|\mbox{\emph{LF}}(N)| = f_L(N) + o(N)$ (the $o$-term may depend on $L$).
\end{theorem}
\proof Consider the family $\mathcal{F}'$ constructed in Proposition \ref{prop8.2}. For each $F' \in \mathcal{F}'$, simply count the sets $A \subseteq [N]$ which satisfy $|A \setminus F'| = o(N)$. Since $|F'| \leq f_L(N)$, the number of such sets $A$ is no more than $2^{f_L(N) + o(N)}$. But \textit{all} sets $A \in \mbox{LF}(N)$ arise from some $F' \in \mathcal{F}'$ in this way and so, since $|\mathcal{F}'| = 2^{o(N)}$, we do indeed have the estimate $|\mbox{LF}(N)| = 2^{f_L(N) + o(N)}$.\endproof\\[11pt]
Observe that the result is best possible apart from the $o$-term, since $\mbox{LF}(N)$ certainly contains all subsets of a strongly $L$-free set with maximal cardinality.\\[11pt]
It turns out that Proposition \ref{prop8.1} can also be derived from Theorem \ref{thm5.3}, though with much weaker quantitative information than that obtainable using the methods of \cite{Gre2,GreRuz1,GreRuz2}. This leads to a unified treatment of the enumeration of $L$-free sets. We sketch the argument here, restricting attention to sum-free sets for simplicity. \\[11pt]
To prove Proposition \ref{prop8.1} it suffices to prove an analagous result with $[N]$ replaced by $G = \mathbb{Z}/p\mathbb{Z}$, where $p \in [2N,4N]$ is a prime. Subsets of $[N]$ may be regarded as subsets of $G$ by reducing mod $p$, and a set $A \subseteq [N]$ is sum-free if and only if $A\pmod{p}$ is sum-free. Let $\epsilon > 0$. For each sum-free set $A \subseteq G$, find an $\epsilon$-regular pair $(R,\eta)$ for $A$. There is some $\rho \in [2\epsilon^{1/3},4\epsilon^{1/3}]$ such that the cardinality of the set
\[ S_{\rho} \; = \; \left\{ x : \rho \leq \alpha_1(x) < \rho + \epsilon^{1/6}\right\} \] is at most $4\epsilon^{1/6} N$. Pick such a $\rho$,
and consider the set $\overline{A}$ consisting of all regular values $x$ such that $\alpha_1(x) \geq \rho$ (note that this differs from the notion of reduced set in Definition \ref{def7.1} only in that we do not restrict ourselves to values of $x$ which also lie in $A$). Take $\mathcal{F}$ to be the collection of all these sets $\overline{A}$. Now it is easy to see that $|A \setminus \overline{A}|$ is small and, using the counting lemma, that $\overline{A}$ is sum-free. It is rather less easy to see that $|\mathcal{F}| = 2^{o(n)}$, and we only give a very brief sketch of the argument.\vspace{11pt}

\noindent In estimating $|\mathcal{F}|$ up to factors of $2^{o(n)}$ we may ignore any non-regular values of $x$ and simply count sets having the form $A^{\circ} = \{x : \alpha_1(x) \geq \rho\}$. Observe that if $\alpha_1(x) \geq \rho + \epsilon^{1/6}$ then $\alpha_1(x + y) \geq \rho$ for all $y \in B_{R,\kappa}$, where $\kappa = \eta\epsilon^{1/6}/20$. This is a consequence of Lemma \ref{lem4.5} (v). Thus $A^{\circ}$ is a union of $x$ for which $x + B_{R,\kappa} \subseteq A^{\circ}$ together with at most $|S_{\rho}|$ extra points. Now by a classical argument of Dirichlet (invoking his principle of the pigeons) $B_{R,\kappa}$ contains an arithmetic progression $P$ of length at least $\kappa  N^{1/|R|}$. Roughly speaking this means that $A^{\circ}$ may be written as a union of longish arithmetic progressions together with $S_{\rho}$. It is then easy to establish a bound $|\mathcal{F}| \leq 2^{c(\epsilon)n}$, where $c(\epsilon) \rightarrow 0$ as $\epsilon \rightarrow 0$.\\[11pt]
\noindent\textit{III. Finding a single regular value: Bourgain's bound for $r_3(n)$.} We owe a large debt to Jean Bourgain and his proof \cite{Bou} of the bound $r_3(n) = O((\log\log n/\log n)^{1/2})$, which is the best currently known. Bourgain's argument essentially amounts, in the language of the present paper, to finding a \textit{single} value of $x$ and a pair $(R,\eta)$ for which $x$ is $\epsilon$-regular. If one is interested in such a weakening of Theorem \ref{thm5.3} then substantial improvements can be made in the bounds. By far the most important difference between our argument and that of Bourgain is that he uses an $\ell^{\infty}$ notion of index in place of our $\ell^2$ definition \eqref{eq11}.\\[11pt]
The fact that a large regular pair can be found was first observed in the context of graph regularity by Koml\'os (unpublished) and elaborated upon by Peng, R\"odl and Rucinski\cite{PRR}. In the arithmetic setting the argument goes through particularly cleanly when $G = (\mathbb{Z}/3\mathbb{Z})^n$, and is then essentially the argument used by Meshulam \cite{Mes}. This was based on the original argument employed by Roth \cite{Rot} to prove that $r_3(n) = o(n)$, but for general groups $G$ Roth's argument does not fit into the framework of regularity since it involves passing to substructures whose size is a small power of $ N$.\\[11pt]
\noindent\textit{IV. Higher arithmetic regularity?} The discussion of III suggests that the analytic proof of $r_3(n) = o(n)$ \cite{Bou,Rot} and the proof via the regularity lemma \cite{FGR,RuzSzem} are perhaps not as different as previously thought. There is also Gowers' analytic proof that $r_4(n) = o(n)$ \cite{Gow2} and a proof via the Frankl-R\"odl regularity lemma for 3-uniform hypergraphs \cite{FraRod}, hereafter termed FRHRL. It would seem to be interesting to ask whether there is an ``arithmetic'' version of FRHRL which is analagous to FRHRL in the same way that Theorem \ref{thm5.3} is analagous to SzRL. A proper understanding of this might be expected to lead to an improvement in the known bounds for $r_4(n)$. Very recently, regularity lemmas for $k$-uniform hypergraphs have been announced independently by Gowers and (various combinations of) Nagle, R\"odl, Schacht and Skokan. Of course, these form part of a more general puzzle.\\[11pt]
The diagram is an attempt to outline what might be hoped for. The three question marks refer to as yet untreated problems. The problem of finding a single regular structure in the hypergraph setting has probably not been investigated, maybe because no potential applications are known. It is quite likely, however, that such a result might not be too difficult to obtain, at least after reading \cite{FraRod} or the more recent works on hypergraph regularity. For that reason it is given the symbol $?^{\ast}$ in our picture.
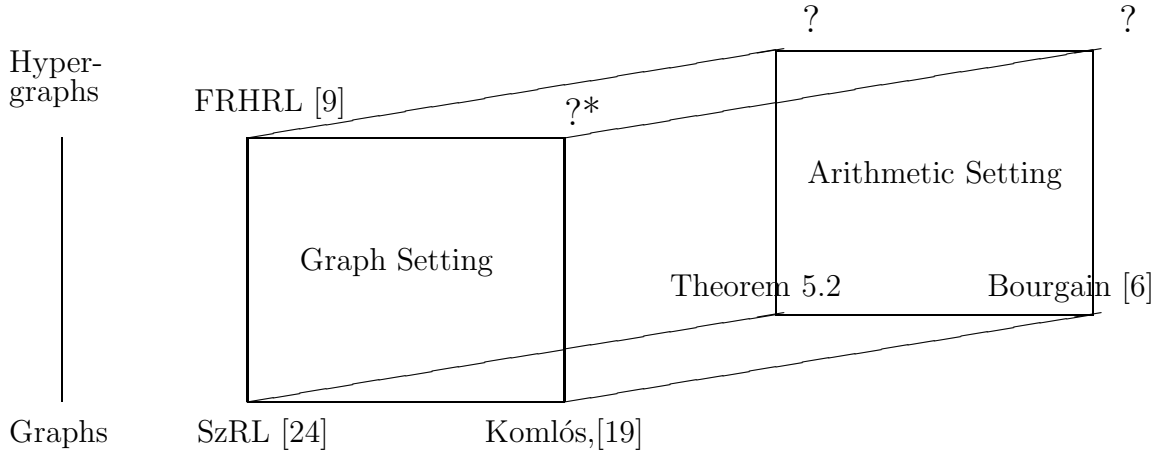
\begin{center}
\begin{picture}(440,250)
\put(90,90){\line(1,0){120}}
\put(90,190){\line(1,0){120}}
\put(90,90){\line(0,1){100}}
\put(210,90){\line(0,1){100}}
\put(90,90){\line(6,1){203}}
\put(210,90){\line(6,1){203}}
\put(90,190){\line(6,1){203}}
\put(210,190){\line(6,1){203}}
\put(290,123){\line(1,0){120}}
\put(290,223){\line(1,0){120}}
\put(290,123){\line(0,1){100}}
\put(410,123){\line(0,1){100}}
\put(110,140){Graph Setting}
\put(302,173){Arithmetic Setting}
\put(71,75){SzRL \cite{Szem}}
\put(180,75){Koml\'os,\cite{PRR}}
\put(70,200){FRHRL \cite{FraRod}}
\put(210,195){\large ?*}
\put(250,130){Theorem \ref{thm5.3}}
\put(370,130){Bourgain \cite{Bou}}
\put(300,230){\large ?}
\put(420,230){\large ?}
\put(20,90){\line(0,1){100}}
\put(0,75){Graphs}
\put(0,215){Hyper-}
\put(0,204){graphs}
\put(90,35){\line(1,0){120}}
\put(0,30){Full regularity}
\put(225,30){One regular object}
\end{picture}
\end{center}
Since the first edition of this paper in October 2003, T. Tao and the author have made some progress on these issues. In particular we have a preprint obtaining a bound $r_4((\mathbb{Z}/5\mathbb{Z})^n) = O(N(\log N)^{-c})$, where $5^n = N$. This is done, essentially, by filling in the top right corner of the above diagram for the particular group $G = (\mathbb{Z}/5\mathbb{Z})^n$. As a result of this work we are able to guess at the correct statements for all of the question marks in the diagram. Proving these in full generality may be a long way off, however.\vspace{11pt}

\noindent One consequence of a suitably formulated ``higher arithmetic regularity lemma'' might be a solution to the following conjecture, which is closely related to Corollary \ref{cor1.6} and the discussion of II. Let $r,m$ be integers with $r \leq m$ and suppose that $L$ is an $r \times m$ matrix of integers with rank $r$. Say that a set $A \subseteq [N]$ is $L$-free if there are no vectors $x \in A^m$ for which $Lx = 0$. 
\begin{conjecture}
Let $A \subseteq [N]$, and suppose that there are $o(N^{m-r})$ vectors in $A^m$ for which $Lx = 0$. Then $A = B \cup C$, where $B$ is $L$-free and $|C| = o(N)$.\end{conjecture}
\noindent Another application of such a result might be a positive answer to Question \ref{BHK-question} for four-term arithmetic progressions.

\section {A Gowers tower for $(\mathbb{Z}/2\mathbb{Z})^n$}\label{sec9} \noindent In this section we show that our regularity results must necessarily give terrible, tower-type bounds. In the context of graphs such a phenomenon was discovered by Gowers \cite{Gow1}: he constructed graphs in which the smallest $\epsilon$-regular partition has a number of parts which grows like $W(\epsilon^{-1/16})$. Our lower bound will not be quite so spectacular but can hardly be described as slowly-growing. Before stating it, we state and prove a preliminary lemma.\\[11pt]
If $M$ is a positive integer, define $F(M)$ to equal $M$ if $M \leq 19$ and $\lfloor M/4\rfloor$ if $M \geq 20$. Define the sequence $\{d_i\}_{i \geq 0}$ by $d_0 = 0$ and $d_{i+1} = F(2^{d_i + d_{i-1} + \dots + d_0})$ for $i \geq 0$. Observe that $d_1 = 1$, $d_2 = 2$, $d_3 = 8$, $d_4 = 512$, $d_5 = 2^{521}$ and that, for $i \geq 6$, $d_i \geq W(i-2)$.
\begin{lemma}\label{lem9.1} Let $M \geq 1$ be an integer, and write $V = (\mathbb{Z}/2\mathbb{Z})^{F(M)}$. Then there are $M$ vectors $\xi_1,\dots,\xi_M \in V$ with the property that any subset of 95 percent of them span $V$.
\end{lemma}
\proof If $M \leq 19$ this is trivial -- simply take $\xi_1,\dots,\xi_M$ to be any basis for $V$. For $M \geq 20$ we use a random approach. Choose the $\xi_i$ independently at random using the uniform distribution on $V$. Let $U$ be a fixed codimension 1 subspace of $V$. The events $\{\xi_i \in U\}$ are independent Bernouilli random variables, and we may invoke a standard tail estimate such as (\cite{AloSpe}, Theorem A.1.4) to deduce that 
\[ \mathbb{P}\left(\mbox{at least 95 percent of the $\xi_i$ lie in $U$}\right) \; \leq \; e^{-M/4}.\]
Thus the probability that some codimension 1 subspace $U$ contains 95 percent of the vectors $\xi_i$ is no more than $2^{F(M)}e^{-M/4}$, which is certainly less than $1$. It follows that there is indeed some choice of the $\xi_i$ satisfying the conclusion of the lemma.\endproof\\[11pt]
Let $G = (\mathbb{Z}/2\mathbb{Z})^n$ and, as usual, write $N = |G|$. The next theorem, which is the main result of this section, provides an example of a function $f : G \rightarrow [0,1]$ such that the largest subgroup $H \leq G$ which is $\epsilon$-regular for $f$ has extremely large index. The definition of what it means for $x$ to be $\epsilon$-regular for $f$ is the obvious one, given what was said in \S \ref{sec2}. A standard probablistic argument such as the one in \cite{Gow1}, Lemma 2 would produce a genuine set $A$ with much the same properties if one was desired.
\begin{theorem} Let $\epsilon < 1/20$ and suppose that $N$ is sufficiently large. Then there is a function $f : G \rightarrow [0,1]$ with the property that any subgroup $H$ which is $\epsilon$-regular for $f$ satisfies $|G/H| \geq W(\textstyle\frac{1}{2}\displaystyle \log_2(1/\epsilon) - 5)$.\end{theorem}
\proof 
Let $s = \lfloor \textstyle\frac{1}{2}\displaystyle \log_2(1/\epsilon) - 2 \rfloor$. Take a fixed nested sequence of subspaces $G  = H_0 \geq H_1 \geq \dots \geq  H_s$ together with subspaces $U_{i+1} \subseteq H_i$
so that $H_i  =  H_{i+1} + U_{i+1}$, where $\dim_{\mathbb{Z}/2\mathbb{Z}} U_i \; = \; d_i$ (the numbers $d_i$ are the ones defined at the start of the section). Define $V_i = U_i + U_{i-1} + \dots + U_0$, so that $H_i + V_i = G$. For each $i$, construct a set $B_i \subseteq G$ as follows. Set $M = 2^{d_i + d_{i-1} + \dots + d_0} = |V_i|$ in Lemma \ref{lem9.1}, and take vectors $\xi_{v} \in U_{i+1}$ indexed by $v \in V_i$, with the property that any $(1 - \epsilon)|V_i|$ of them span $U_{i+1}$. For each $v$, let $\xi_v^{\circ} = \{ u \in U_{i+1} : \langle u,\xi_v \rangle = 0\}$.
Writing $G = H_{i+1} + U_{i+1} + V_i$, let $B_i$ be the set defined by
\[ B_i \; = \;  \bigcup_{v \in V_i} \left(H_{i+1} + \xi_v^{\circ} + v\right).\]
$B_i$ has cardinality $N/2$, and consists of a codimension one subspace in each coset of $H_i$. Define \[ f \; = \; \textstyle\frac{1}{2}\displaystyle \bigg(B_0 + \frac{1}{4}B_1 + \dots + \frac{1}{4^s}B_s\bigg).\]
We will prove, by induction on $i$, that if $H$ is $\epsilon$-regular for $f$ then $H \subseteq H_i$. Suppose we know that $H \subseteq H_i$. Take $g \in G$ and write it as $y + v$, where $y \in H_i$ and $v \in V_i$. Consider the set $(B_i)^{+g}_H$. We have
\[ (B_i)^{+g}_H(x) \; = \; \left\{\begin{array}{ll} 1 & \mbox{if $x \in H \cap \left(H_{i+1} + \xi_{v}^{\circ} + y\right)$} \\ 0 & \mbox{otherwise}. \end{array}\right.\]
Now $\xi_v^{\circ} + H_{i+1} + y$ is a hyperplane (coset of a codimension 1 subspace) in $H_i$ which is perpendicular to $\xi_v$. Therefore its intersection with $H$ is either empty, all of $H$ or a hyperplane in $H$. The latter case arises when $H$ is not a subspace of $H_{i+1} + \xi_v^{\circ}$. Supposing this is the case, then $|\widehat{(B_i)^{+g}_H}(\xi_v)|  =  |H|/2$.
Furthermore if $j < i$ then $B_j$ is, by construction, a union of $H_i$-cosets. Since $H \subseteq H_i$ it follows that $(B_j)^{+g}_H$ is either empty or else all of $H$, and so $\widehat{(B_j)^{+g}_H}(\xi_v) = 0$. When $j > i$ we have the trivial bound $|\widehat{(B_j)^{+g}_H}(\xi_v)| \leq |H|$. It follows that if $H \not \subseteq H_{i+1} + \xi_v^{\circ}$ then
\[ |H|^{-1}|\widehat{f}^H_g(\xi_v)| \; \geq \; \textstyle\frac{1}{4}\displaystyle\cdot 4^{-i}  - \textstyle\frac{1}{2}\displaystyle \sum_{j > i} 4^{-j} \; \geq \; \textstyle\frac{1}{16}\displaystyle 4^{-i} \; \geq \; \textstyle\frac{1}{16}\displaystyle 4^{-s} \;  > \; \epsilon.\]
As $H$ is assumed to be $\epsilon$-regular for $f$ this can hold for no more than $\epsilon N$ values of $g \in G$. This means that the proportion of $v \in V_i$ for which $H \not\subseteq H_{i+1} + \xi_v^{\circ}$ is at most $\epsilon$. Let $\Omega$ be the remaining values of $v$, of which there are at least $(1 - \epsilon)|V_i|$. By the choice of the vectors $\xi_v$, the collection $\{\xi_v\}_{v \in \Omega}$ spans $U_i$. Thus $H$, which lies in
\[ \bigcap_{v \in \Omega} \left(H_{i+1} + \xi_v^{\circ}\right),\] must in fact be contained in $H_{i+1}$.\\[11pt]
We have completed the inductive step. It follows that if $H$ is $\epsilon$-regular for $f$ then $H \subseteq H_s$, which leads to the lower bound on $|G/H|$ stated in the theorem.\endproof

\section{Acknowledgements} \noindent The author would like to thank Tim Gowers and Imre Ruzsa for unwittingly making remarks that helped the author have the ideas in this paper, and to the journal for encouraging him to find a further application (Theorem \ref{thm1.9}) of the regularity lemma.

\providecommand{\bysame}{\leavevmode\hbox to3em{\hrulefill}\thinspace}
\providecommand{\MR}{\relax\ifhmode\unskip\space\fi MR }
\providecommand{\MRhref}[2]{%
  \href{http://www.ams.org/mathscinet-getitem?mr=#1}{#2}
}
\providecommand{\href}[2]{#2}

     \end{document}